\documentclass[11pt]{article}
\usepackage{latexsym}
\usepackage{amsmath}
\usepackage{cite}
\usepackage{amsthm,color,bm}
\usepackage{amssymb}
\usepackage{mathrsfs}
\usepackage{indentfirst}
\usepackage{times,cases}
\usepackage{dsfont}

\usepackage{lscape}
\newtheorem{theorem}{\noindent Theorem}[section]
\newtheorem{proposition}[theorem]{\noindent Proposition}
\newtheorem{definition}[theorem]{\noindent Definition}
\newtheorem{lemma}[theorem]{\noindent Lemma}
\newtheorem{remark}[theorem]{\noindent Remark}
\newtheorem{corollary}[theorem]{\noindent Corollary}
\numberwithin{figure}{section}
\numberwithin{equation}{section}

\renewcommand{\theequation}{\thesection.\arabic{equation}}

\newcommand{\cB}{{\mathcal B}}
\newcommand{\cC}{{\mathcal C}}
\newcommand{\cD}{{\mathcal D}}

\newcommand{\cF}{{\mathcal F}}

\newcommand{\cI}{{\mathcal I}}
\newcommand{\cK}{{\mathcal K}}

\newcommand{\cP}{{\mathcal P}}
\newcommand{\cR}{{\mathcal R}}

\newcommand{\cX}{{\mathcal X}}
\newcommand{\cY}{{\mathcal Y}}
\newcommand{\cZ}{{\mathcal Z}}

\newcommand{\sA}{{\mathscr A}}

\newcommand{\sD}{{\mathscr D}}

\def\R{\mathbb{R}}
\def\bE{\mathbb{E}}
\def\bP{\mathbb{P}}

\def\N{\mathbb{N}}
\def\1{\mathds{1}}
\newcommand{\di}{\mathrm{d}}
\newcommand{\me}{\mathrm{e}}

\def\disp{\displaystyle}

\def\bc{\begin{center}}
\def\ec{\end{center}}
\def\be{\begin{equation}}
\def\ee{\end{equation}}
\def\bea{\begin{eqnarray}}
\def\eea{\end{eqnarray}}
\def\ba{\begin{array}}
\def\ea{\end{array}}
\def\benu{\begin{enumerate}}
\def\eenu{\end{enumerate}}
\def\bt{\begin{theorem}}
\def\et{\end{theorem}}
\def\bl{\begin{lemma}}
\def\el{\end{lemma}}
\def\bco{\begin{corollary}}
\def\eco{\end{corollary}}
\def\bn{\begin{numcases}}
\def\en{\end{numcases}}
\def\br{\begin{remark}}
\def\er{\end{remark}}
\def\bd{\begin{definition}}
\def\ed{\end{definition}}
\def\bp{\begin{proposition}}
\def\ep{\end{proposition}}
\def\bo{\begin{proof}}
\def\eo{\end{proof}}
\def\bx{\begin{example}}
\def\ex{\end{example}}
\def\bal{\begin{align}}
\def\eal{\end{align}}

\def\pa{\partial}
\def\al{\alpha}\def\b{\beta}
\def\De{\Delta} \def\de{\delta}
\def\na{\nabla}
\def\lam{\lambda} 

\def\ve{\varepsilon}
\def\sig{\sigma}
\def\vsig{\varsigma}
\def\vp{\varphi}
\def\w{\omega}\def\W{\Omega}
\def\gam{\gamma}\def\Gam{\Gamma}

\def\~{\widetilde}

\def\ol{\overline}

\def\Cap{\bigcap}
\def\Cup{\bigcup}

\def\ra{\rightarrow}

\def\8{\infty}
\def\X{\times}

\def\mb{\mbox}
\def\di{{\rm d}}
\def\me{{\rm e}}

\def\suo{\!\!\!}
\def\es{\emptyset}

\def\sm{\setminus}

\def\hs{\hspace{0.4cm}}
\def\Vs{\vskip10pt}
\def\vs{\vskip5pt}

\def\({\left(}
\def\){\right)}

\begin{document}


\begin{center}
    {\large \bf Upper semi-continuity of random attractors and existence of invariant measures
    for nonlocal stochastic Swift-Hohenberg equation with multiplicative noise}
\vspace{0.5cm}\\
{Jintao Wang$^{1}$,\quad Chunqiu Li$^{2,*}$,\quad Lu Yang$^{3}$,\quad Mo Jia$^{4}$}\\\vspace{0.3cm}

{\small $^{1,2}$Department of Mathematics, Wenzhou University, Wenzhou 325035, China\\\vspace{0.2cm}
        $^{3}$School of Mathematics and Statistics, Lanzhou University, Lanzhou, 730000, China\\\vspace{0.2cm}
        $^{4}$School of Mathematical Science, Tianjin University, Tianjin, 300350, China}
\end{center}


\renewcommand{\theequation}{\arabic{section}.\arabic{equation}}
\numberwithin{equation}{section}


\begin{abstract}
In this paper, we mainly study the long-time dynamical behaviors of 2D nonlocal stochastic Swift-Hohenberg equations with multiplicative noise from two perspectives.
Firstly, by adopting the analytic semigroup theory, we prove the upper semi-continuity of random attractors in the Sobolev space $H_0^2(U)$, as the coefficient of the multiplicative noise approaches zero.
Then, we extend the classical ``stochastic Gronwall's lemma", making it more convenient in applications.
Based on this improvement, we are allowed to use the analytic semigroup theory to establish the existence of ergodic invariant measures.
\vs

\noindent\textbf{Keywords:} Nonlocal stochastic Swift-Hohenberg model; multiplicative
noise; random attractor; asymptotic compactness; upper semi-continuity; Markovian transition semigroup; Feller property; ergodic invariant measures.
\vs

\noindent{\bf AMS Subject Classification 2010:}\, 60H15, 86A05, 37L55, 37A25

\end{abstract}

\vspace{-1 cm}

\footnote[0]{\hspace*{-7.4mm}
$^{*}$ Corresponding author.\\
E-mail address: wangjt@hust.edu.cn (J.T. Wang); lichunqiu@wzu.edu.cn (C.Q. Li);
yanglu@lzu.edu.cn (L. Yang);\\ jiamomath@tju.edu.cn (M. Jia).}


\section{Introduction}

In the research of geophysical fluid flows in atmosphere, oceans and the earth's mantle,
there often arise fluid convective phenomena resulting from the density variations.
A prototypical model to describe these phenomena is the Rayleigh-B\'enard convection,
to forecast the spatio-temporal convection patterns.
The mathematical model for the Rayleigh-B\'enard convection involves nonlinear Navier-Stokes
partial differential equations coupled with the temperature equation.
When the Rayleigh number approaches the onset of the convection, the Rayleigh-B\'enard convection model
may be approximately reduced to an amplitude or order parameter equation as derived by
Swift and Hohenberg (\cite{SH77}), after whom the Swift-Hohenberg equation (SHE for short) was named.
The SHE has featured in different branches of physics, ranging from hydrodynamics to nonlinear optics,
such as the Taylor-Couette flow (\cite{PM80}), plasma confinement in toroidal devices (\cite{HS92,LMR75}),
viscous film flow, lasers(\cite{LMN94}) and pattern formation.

The localized one-dimensional version (see \cite{SH77}) of the SHE reads as follows,
$$
u_t=au-(1+\pa_{xx})u-u^3,
$$
where cubic term $u^3$ is used as an approximation of a nonlocal integral term.
When reexamining the rationale for using the Swift-Hohenberg model as a
reliable model of the spatial pattern evolution in specific physical systems,
Roberts in \cite{R92,R95} argued that, although the localization approximation used in \cite{SH77}
makes some sense in the one-dimensional case, this approximation is deficient for the two-dimensional one.
Therefore one should use the nonlocal SHE (\cite{R92,R95,SH77}),
$$u_t=au-(1+\De)^2u+uG*u^2,$$
where
\be\label{1.0}[G*u^2](x_1,x_2,t)=\int_U G\(\sqrt{(x_1-y_1)^2+(x_2-y_2)^2}\)u^2(y_1,y_2,t)\di y_1\di y_2,\ee
$u(x_1,x_2,t)$ is the unknown amplitude function, $a$ measures the difference of the Rayleigh number
from its critical onset value, $\De=\pa_{x_1x_1}+\pa_{x_2x_2}$ is the Laplace operator,
$G:\R^+\ra\R$ is radially symmetric on $\R^2$ and $*$ denotes the convolution.
The equation is defined for $t>0$ and $(x_1,x_2)\in U$, where $U$ is a bounded planar domain.

In spite of the special importance of the nonlocal case in physical facts,
there have been few investigations (\cite{G15,LGDE00,WSD05}) on the nonlocal SHE up to now,
comparing the abundant works having been done on other types of SHE recently;
see \cite{LGDE00,WSD05,CHHL17,GDL16,LWZ20,S18,WYD20,XG10}
and the references therein.

As another physical fact, fluid systems are often affected by random environmental influences.
Due to the delicate impact on the deterministic cases, increasing concern of the randomness has occurred
all over the development of evolution systems
(e.g. \cite{G15,WSD05,EH01,EKZ17,GDL16,JQW18,GY19,DLS03,GKV14,DaPZ14,MSY16,GHZ09,LWZ20,XHZ19,XZH19,CLR98,CKW18,LX20})
in recent decades.
It becomes centrally significant to take stochastic effects into account for mathematical models of
complex phenomena in engineering and science.
The involvement of random motions makes the theoretical results more reasonable and reliable in the sense of statistics.
Stochastic partial differential equations are appropriate models for randomly influenced spatially extended systems.
Furthermore, the participation of such effects has brought interesting new mathematical problems at the interface of
probability and partial differential equations.
\Vs

In this paper, we consider the limiting behaviors of random attractors and invariant measures of
the following planar nonlocal stochastic SHE with multiplicative noise,
\bn{}
\di u+[\De^2u+2\De u+au+uG*u^2+h(x)]\di t=\epsilon u\di W_t,\;\;x\in U,~ t>0,\label{1.1}\\
u=\frac{\pa u}{\pa\bm n}=0,\;\; x \in \pa U,~t\geqslant0,\label{1.2}\\
u(0,x)=u_0(x),\;\;x\in U,\label{1.3}
\en
where $a\in\R$, $h$ is the external forcing term with $h\in L^2(U)$, $G*u^2$ is given in \eqref{1.0},
$\epsilon$ is a small positive parameter, $\bm n$ is the external normal vector on the boundary $\pa U$,
$u_0$ is the initial datum and $W_t$ is a two-sided real-valued Wiener process
on the probability space to be given later.
As mentioned above, we have only discovered three references \cite{G15,LGDE00,WSD05} that research on nonlocal SHEs.
Lin et al in \cite{LGDE00} compared the difference of bounds for the dimensions of the global attractors for
local and nolocal two-dimensional autonomous SHEs.
Wang, Sun and Duan studied in \cite{WSD05} the local and nonlocal SHEs under the influence of additive white noise.
Recently, the existence of random attractors and invariant manifolds for a nonlocal stochastic SHE was
considered by Guo in \cite{G15}.

We compose our work mainly into two topics under two distinct assumptions on the non-negative kernel $G$ ---
the positive kernel and a special non-negative kernel, which will be presented in Section 2 in details.

Our first concern is the limiting behaviors of random attractors of the nonlocal
stochastic SHE \eqref{1.1} -- \eqref{1.3}, when the coefficient of stochastic term tends to 0.
The upper semi-continuity of random attractors illustrates the decreasing of effective influences of randomness
on the long-time dynamical behaviors as the intensity of random action imposed on the system fades away.
This type of continuity for attractors has been studied by many authors;
see e.g. \cite{CKW18,CLR98,HKU20,LX20} and the references therein.
The existence and regularity of random attractors of the problem \eqref{1.1} -- \eqref{1.3} can be obtained by a similar routine in \cite{G15}.

Our aim in this part is to show the upper semi-continuity of random attractors in $H_0^2(U)$.
Here what matters is the asymptotic compactness of the generated process in $H_0^2(U)$,
since the boundary condition \eqref{1.2} makes it difficult to treat the term $\int_U\De^mu\De^2u\di x$
for each $m\geqslant2$ through integration by parts.
Thanks to the analytic semigroup theory of fractional power of the infinitesimal generator,
saying, the operator $\De^2$ on $H_0^2(U)\cap H^4(U)$, we eventually overcome this technical problem.

The second issue of this paper is the existence of ergodic invariant measures.
The invariant measure is an important tool to study long-time behaviors of the solutions to
stochastic dynamical systems.
There have been plenty of works in this field, including invariant measures for stochastic systems
(see \cite{WSD05,EH01,EKZ17,GKV14,DaPZ14,MSY16}) and deterministic systems
(see \cite{CG12,LR14,XL20,WZC20,ZXL18,ZY17}).
On invariant measures for infinite-dimensional stochastic dynamical systems,
most of the historical discussions concerned additive noise for the stochastic term
(e.g. \cite{WSD05,EH01,EKZ17}) and relative techniques have been well developed.
Early in 2005, Wang, Sun and Duan in \cite{WSD05} studied the existence of ergodic invariant measures
for nonlocal stochastic SHE with additive noise.
The authors gave an exhaustive investigation for the generated semigroup on the metric space of
all Borel probability measures in $L^2(U)$.
They proved that when sufficiently many of its Fourier modes are forced,
the system has a unique invariant measure, or equivalently, the dynamics is ergodic.
However, it is more challenging for the case of multiplicative noise,
since the stochastic term itself also changes with the evolution of the unknown function,
whose estimate requires different tricks.

In our situation, we consider the nonlocal stochastic SHE with an external force $h$ and multiplicative noise.
We follow the classical Krylov-Bogoliubov procedure to show the existence of invariant measures on $L^2(U)$,
and use \cite[Proposition 11.12]{DaPZ14} and the Krein-Milman Theorem (\cite{Conway90}) to prove the existence of ergodic invariant measures.

Before this procedure, we need to confirm that the corresponding Markovian transition semigroup $\cP_t$ is Feller.
This can be essentially solved by the well-known ``stochastic Gronwall's lemma", Proposition \ref{p4.1}.
However, since $h\in L^2(U)$, when we take the expectation of $H_0^2(U)$-norm,
it seems difficult to apply $-\De$ to \eqref{1.1} directly
and this prevents us making good use of the It\^o's formula.
We are compelled to choose the integral equation (see \eqref{4.A} below) and the analytic semigroup theory to surmount this obstacle.
In this way, the classical ``stochastic Gronwall's lemma" can not directly apply to our estimation.
We hence develop the ``stochastic Gronwall's lemma" in Lemma \ref{le4.2}, which enables us to handle this difficulty.
Actually, instead of the arbitrariness of the stopping times $\tau'$ and $\tau''$ in the condition \eqref{4.1} in Proposition \ref{p4.1},
Lemma \ref{le4.2} only requests the inequality \eqref{4.1} to hold when the interval between $\tau'$ and $\tau''$ is sufficiently short.
This is evidently much easier to check.
Moreover, the relaxed condition \eqref{4.2} is turned out to be equivalent to the previous one \eqref{4.1}.

The remainder of this paper is arranged as follows.
In the subsequent section, we set some basic notations and then recall the existence and uniqueness of the solutions to \eqref{1.1} -- \eqref{1.3}
under two kinds of kernels --- the positive kernel and a special non-negative kernel.
The third section is devoted to the verification of the upper semi-continuity of the random attractors in $H_0^2(U)$ for the nonlocal stochastic SHE with multiplicative noise.
We develop the ``stochastic Gronwall's lemma" and give the existence of ergodic invariant measures in $L^2(U)$ in Section 4.
In the last section, we summarize this paper and propound some related interesting mathematical problems.

\section{Existence and Uniqueness of Solutions}\label{s2}

We first set some basic notations and properties that will be used frequently in this paper.
\subsection{Basic notations and properties}

For metric spaces $X$ and $Y$, we conventionally denote by $\cC(X,Y)$ ($\cC_{\rm b}(X,Y)$)
the collection of continuous (and bounded) functionals from $X$ to $Y$.
When $Y=\R$, we simply use $\cC(X)$ ($\cC_{\rm b}(X)$) to represent $\cC(X,\R)$ ($\cC_{\rm b}(X,\R)$).

Denote the scalar product of the Hilbert space $L^2(U)$ by
$$(u,v)=\int_U uv\di x,$$
and the norm simply by $\|\cdot\|$.
For the space $H_0^2(U)$ and the negative Laplacian $-\De$ on $H_0^2(U)$, it is well known that,
the norm of $H^2_0(U)$ defined as
$$\|u\|_{H^2_0(U)}=\|\De u\|,$$
is equivalent to the usual norm of $H_0^2(U)$.
We write $H=L^2(U)$ and $V=H_0^2(U)$ for simplicity.
\Vs

For sake of computational convenience afterwards, we recall some basic knowledge of the fractional power of
sectorial operator (see \cite{Hen81,LW18,LZ98}) as follows.

Definitely, the operator $A:=\De^2$ is a sectorial operator in $H$ with its domain
$$\sD(A)=V\cap H^4(U).$$
The operator $A$ is positive, self-adjoint and owns a basis of eigenfunctions $\{w_i\}_{i\in\N_+}$
which is orthonormal in $H$ and associated with the eigenvalues $\{\lam_i\}_{i\in\N_+}$ such that
$$0<\lam_1<\lam_2\leqslant\lam_3\leqslant\cdots\leqslant\lam_i\leqslant\ra+\8.$$
For all $\mu\in\R$, we define $A^\mu$ by setting
$$
A^\mu\phi=\sum_{i=1}^\8\lam_i^\mu\xi_iw_i,\hs\mb{when }\phi=\sum_{i=1}^\8\xi_iw_i\in\sD(A^\mu)\hs
\mb{with}\hs\sD(A^\mu)=\{\phi\in H:A^\mu\phi\in H\}.
$$
The space $\sD(A^\mu)$ is thus a Banach space with the norm $\|\cdot\|_\mu$ such that
$$\|\phi\|_\mu=\|A^{\mu}\phi\|.$$

It is well known that $-\De=A^{\frac12}$, $V=\sD(A^{\frac12})$ and
the embedding $\sD(A^{\mu_1})\subset\sD(A^{\mu_2})$ is compact for all $\mu_1>\mu_2$.
If $\phi\in\sD(A^\mu)$, then
\be\label{2.1}A^\mu\me^{-At}\phi=\me^{-At}A^\mu\phi\hs\mb{for all }t\geqslant0.\ee
Let $\lam_0\in(0,\lam_1)$.
Then there are positive constants $C$ and $C_\mu$ for all $\mu>0$, such that
\be\label{2.2}\|\me^{-At}\|\leqslant C\me^{-\lam_0t},\hs
\|A^\mu\me^{-At}\|\leqslant C_\mu t^{-\mu}\me^{-\lam_0t}\hs\mb{for all }t>0.\ee

\subsection{Existence of global solutions}

Now we consider the existence and uniqueness of the solutions for the two-dimensional nonlocal stochastic
Swift-Hohenberg model \eqref{1.1} -- \eqref{1.3} with multiplicative noise.
In order to do this, we will somehow transform the original stochastic problem into a cocycle driven
by an ergodic metric dynamical system.

In the sequel, we consider the probability space $(\W,\cF,\bP)$, where
$$
\W=\left\{\w\in \cC(\R,\R) : \w(0)=0 \right\},
$$
$\cF$ is Borel $\sig$-algebra induced by the compact-open topology of $\W$ and $\bP$ is
the corresponding Wiener measure on $(\W,\cF)$ (see \cite{A98}).
Define the time shift by
$$
\theta_t\w(\cdot)=\w(\cdot+t)-\w(t), \quad\w\in\W,\,t\in\R.
$$
Then $(\W,\cF,\bP,(\theta_t)_{t\in\R})$ is an ergodic metric dynamical system.
Note that, there exists a $\theta_t$-invariant set $\~\W\subseteq\W$ of full $\bP$ measure such that for each $\w\in\~\W$,
$$
\frac{\w(t)}{t}\ra 0\quad\mb{as}~t\ra\pm\8.
$$
In the following, we only consider the space $\~\W$ instead of $\W$, and write $\~\W$ as $\W$ for convenience.

We identify $W_t(\w)=\w(t)$ with $\w\in\W$ defined as above.
Consider the stochastic stationary solution to the one-dimensional Ornstein-Uhlenbeck equation
\be\label{2.3}
\di z+z\di t=\di W_t.
\ee
By \cite{DLS03,LWZ20}, we know that the solution
\be\label{2.4}t\mapsto z(\theta_t\w):=-\int_{-\8}^0\me^s(\theta_t\w)(s)\di s,\hs t\in\R\ee
to \eqref{2.3} is continuous in $t$ for every $\w\in\W$ and that the random variable $|z(\theta_t\w)|$ is tempered, satisfying
\be\label{2.5}\lim_{t\ra\pm\8}\frac{|z(\theta_t\w)|}{t}=0,\hs\lim_{t\ra\pm\8}\frac1t\int_0^tz(\theta_s\w)\di s=0.\ee
Moreover, from \cite[Proposition 4.3.3]{A98}, it follows that there exists a tempered variable $r(\w)>0$ such that, for $\bP$-a.s. $\w\in\W$,
$$|z(\theta_t\w)|\leqslant\me^{\frac{|t|}{2}}r(\w),\hs t\in\R.$$
\vs

Now we set
\be\label{2.6}v(t)=\me^{-\epsilon z(\theta_t\w)}u(t),\ee
where $u$ is a solution to the problem \eqref{1.1} -- \eqref{1.3} and $z(\theta_t\w)$ is given by \eqref{2.4}.
Then $v$ satisfies the following problem,
\bn{}
\frac{\pa v}{\pa t}+\De^2v+2\De v+av+\me^{2\epsilon z(\theta_t\w)}vG*v^2+\me^{-\epsilon z(\theta_t\w)}h=\epsilon z(\theta_t\w)v,\;\;x\in U,~ t>0,\label{2.7}\\
v=\frac{\pa v}{\pa\bm n}=0,\;\; x \in \pa U,~t\geqslant0,\\
v(0,x)=v_0(x)=\me^{-\epsilon z(\w)}u_0(x),\;\;x\in U.\label{2.9}
\en

Basing on different choices of the kernel $G$, we will consider the nonlocal stochastic SHEs with the kernel $G$ under two different conditions in the following, namely, positive kernel and a special non-negative kernel.
\vs

\noindent{\bf\large The case of positive kernel}

For the nonlocal stochastic SHE \eqref{1.1} with positive kernel $G$,
we assume that\\
{\bf (Gp)} the functions $G,\,|\na G|,\,\De G\in L^\8(U)$ such that there are positive constants $\al$ and $\b$ such that for every $x\in U$,
\be\label{2.10}G(|x|)\geqslant\al\hs\mb{and}\hs\|G\|_\8,\;\|\na G\|_\8,\;\|\De G\|_\8\leqslant\b,\ee
where $\na=(\pa_{x_1},\pa_{x_2})$ is the gradient operator.
\vs

\noindent{\bf\large The case of non-negative kernel}

For the nonlocal stochastic SHE \eqref{1.1} with non-negative kernel $G$,
we consider a special kernel, which is a mollifier.
Define
$$J(r)=\left\{\ba{ll}c\exp\(-\frac1{1-r^2}\),&r<1,\\0,r\geqslant1,\ea\right.\hs\mb{where}\hs c=\(\int_{D_1}\exp\(-\frac1{1-r^2}\)\di x\)^{-1},$$
$r=\sqrt{x_1^2+x_2^2}$ and $D_1=\{(x_1,x_2)\in\R^2:r<1\}$.
Define also for some $\varrho>0$, $J_\varrho(r)=\varrho^{-2}J(r/\varrho)$.
Let $\cC_0(\ol U)=\{\psi\in\cC(\ol U):{\rm supp}\psi\subset U\}$.
For every $\psi\in \cC_0(\ol U)$, define the mollifier of $\psi$ by the convolution $J_\varrho*\psi$.
It is well known that
$$\|J_\varrho*\psi-\psi\|_{\cC_0(\ol U)}\ra0\hs\mb{as }\varrho\ra0.$$
Therefore, for each given $\ve>0$, there is a positive constant $\varrho_0$ such that
$$J_\varrho*\psi\geqslant\psi-\ve\hs\mb{as }\varrho\in(0,\varrho_0).$$
Moreover, $J_{\varrho_0}$ satisfies
$$0\leqslant J_{\varrho_0}\leqslant\frac{c}{\varrho_0^2},\hs\mb{and}\hs |\na J_{\varrho_0}|,\,\De J_{\varrho_0}\in L^\8(U).$$
By considering the special kernel $G$ to be $J_{\varrho_0}$ in \eqref{1.1} with $\varrho_0$ given above,
we can endow $G$ the following assumption that\\
{\bf (Gn)} there are positive constants $\al$, $\b$ and $\de$ such that for every $x\in U$, $G$ satisfies
$$G*\psi\geqslant\al \psi-\de\mb{ for all }\psi\in\cC_0(\ol U)\;\mb{ and}\hs \|G\|_\8, \|\na G\|_\8,\;\|\De G\|_\8\leqslant\b.$$
\Vs

Under the assumption either {\bf(Gp)} or {\bf(Gn)} and following standard arguments (see \cite{DaPZ14,Hen81,Tem88}),
for every initial datum $v_0(\w)\in H$, the system \eqref{2.7}-\eqref{2.9} possesses a unique global strong solution (\cite{G15,LGDE00})
\be\label{2.11} v(\cdot;\w,v_0(\w))\in\cC([0,+\8);H)\cap L^\8(0,+\8;V)\cap L^2(0,T;\sD(A)),\ee
which is continuous with respect to $v_0(\w)\in V$ for all $t\geqslant0$.
In particular, under the assumption {\bf(Gn)} and following the standard regularity argument in \cite{Hen81} with the Sobolev embedding theorem,
we have moreover the solution $u(t)\in \cC_0(\ol U)$ for all $t\geqslant0$.
Then we can deduce from {\bf(Gn)} that
\be\label{2.12}G*u^2\geqslant\al u^2-\de\hs\mb{and}\hs \|G\|_\8, \|\na G\|_\8,\;\|\De G\|_\8\leqslant\b.\ee

Let
\be\label{2.13}u(t;\w,u_0)=\me^{\epsilon z(\theta_t\w)}v(t;\w,\me^{-\epsilon z(\w)}u_0).\ee
Then according to the result in \cite[Lemma 2.2]{DLS03}, for each $u_0\in H$, \eqref{2.13} gives the problem \eqref{1.1} -- \eqref{1.3}
a unique (in almost-sure sense) global {\em strong solution}, such that
$$u(t;\w,u_0)+\int_0^t[Au(s)+2\De u(s)+au(s)+u(s)(G*u^2)(s)+h(x)]\di s=u_0+\int_0^t\epsilon u(s)\di W_s$$
makes sense.
Moreover, the solution $u(t,\w,\cdot)$ is continuous as a mapping from $V$ to $V$ for all $t\geqslant0$.

In the following estimates, we denote $c$ as an arbitrary positive constant, which only depends on the parameters of
the original problem and the assumptions, i.e., $a,\al,\beta,U$, $\de$ (in case when {\bf(Gn)} holds for $G$),
$\mu$ (used in Section 3) and $p$ (appearing in Section 4), and may be different from line to line and even in the same line.
We also use $a\lesssim b$ to denote the $a\leqslant cb$ for simplicity, for the $c$ introduced above.

\section{Upper Semi-continuity of Random Attractors}

We first recall basic concepts and results related to random attractors for random dynamical systems; see \cite{A98,CLR98} for details.

\subsection{Preliminaries on random attractors}

Let $(\W,\cF,\bP,(\theta_t)_{t\in\R})$ be a measurable dynamical system, which can be viewed as that one we have introduced in the last section.
Let $(X,\|\cdot\|_X)$ be a Banach space with Borel $\sig$-algebra $\cB(X)$.
\bd A {\bf continuous random dynamical system} on $X$ over $(\W,\cF,\bP,(\theta_t)_{t\in\R})$ is a measurable mapping
$$\Phi:\R^+\X\W\X X\ra X,\hs(t,\w,x)\mapsto\Phi(t,\w,x)$$
such that for $\bP$-a.s. $\w\in\W$,
\benu\item[(i)] $\Phi(0,\w,\cdot)={\rm id}$ on $X$,
\item[(ii)] $\Phi(t+s,\w,\cdot)=\Phi(t,\theta_s\w,\Phi(s,\w,\cdot))$ for all $s,\,t\in\R^+$ and $\w\in\W$, and
\item[(iii)] $\Phi(t,\w,\cdot):X\ra X$ is continuous.
\eenu\ed

A {\em random set} $\widehat D=\{D(\w)\}_{\w\in\W}$ is a family of subsets of $X$ indexed by $\w$ such that for every $x\in X$ the mapping $\w\mapsto \di(x,D(\w))$ is measurable with respect to $\cF$.
Let $\widehat D'=\{D'(\w)\}_{\w\in\W}$.
We write ``$\widehat D'\subseteq \widehat D$'' to indicate that $D'(\w)\subseteq D(\w)$ for each $\w\in\W$.

A {\em random compact set} $\widehat D=\{D(\w)\}_{\w\in\W}$ is a random set such that $D(\w)$ is compact for all $\w\in\W$.
A random set $\widehat D$ is said to be {\em bounded} if there exists a $u_0\in X$ and a random variable $R(\w)>0$ such that
$$D(\w)\subset\{u\in X:\|u-u_0\|_X\leqslant R(\w)\}\hs\mb{ for all }\w\in\W.$$
A random set $\widehat D=\{D(\w)\}_{\w\in\W}$ is said to be {\em invariant}, if $\Phi(t,\w,D(\w))=D(\theta_{t}\w)$ for all $t\geqslant0$.
A random variable $\{b(\w)\}_{\w\in\W}$ is called {\em tempered} with respect to $(\theta_t)_{t\in\R}$ if for $\bP$-a.s. $\w\in\W$,
$$\lim_{t\ra\8}\me^{-ct}b(\theta_{-t}\w)=0\hs\mb{for all }c>0.$$
A bounded random set $\{B(\w)\}_{\w\in\W}$ of $X$ is called {\em tempered} with respect to $(\theta_t)_{t\in\R}$ if $\|B(\w)\|_X$ is tempered with respect to $(\theta_t)_{t\in\R}$, where $\|B\|_X=\sup_{x\in B}\|x\|_X$.

Let $\cD$ be a collection of random sets of $X$.
We say $\cD$ is {\em inclusion-closed} provided that, for arbitrary random sets $\widehat D$ and $\widehat D'$, the conditions $\widehat D'\subseteq \widehat D$ and $\widehat D\in\cD$ ensure $\widehat D'\in\cD$.

\bd A random dynamical system $\Phi$ is said to be {\bf pullback $\cD$-asymptotically compact} in $X$, provided for $\bP$-a.s. $\w\in\W$, the sequence $\Phi(t_n,\theta_{-t_n}\w,x_n)$ has a convergent subsequence in $X$ whenever $t_n\ra+\8$, and $x_n\in B(\theta_{-t_n}\w)$ with $\{B(\w)\}_{\w\in\W}\in\cD$.\ed

\bd A random set $\widehat K=\{K(\w)\}_{\w\in\W}\in\cD$ is called a {\bf random absorbing set} for $\Phi$ in $\cD$ if for every $\widehat B\in\cD$ with $\widehat B=\{B(\w)\}_{\w\in\W}$ and $\bP$-a.s. $\w\in\W$, there exists $t_{\widehat B}(\w)>0$ such that
$$\Phi(t,\theta_{-t}\w,B(\theta_{-t}\w))\subset K(\w)\hs\mb{for all }t\geqslant t_{\widehat B}(\w).$$\ed
\bd Let $\cD$ be a collection of random sets of $X$.
A random compact invariant set $\widehat\sA$ with $\widehat\sA=\{\sA(\w)\}_{\w\in\W}$ of $X$ is called a {\bf $\cD$-random attractor} for $\Phi$, provided that for $\bP$-a.s. $\w\in\W$, $\sA(\w)$ pullback attracts every random set in $\cD$, that is, for every $\{B(\w)\}_{\w\in\W}\in\cD$,
$$\lim_{t\ra\8}{\rm dist}_X(\Phi(t,\theta_{-t}\w, B(\theta_{-t}\w)),\sA(\w))=0,$$
where ${\rm dist}_X$ is the Hausdorff semi-distance, i.e., ${\rm dist}_X(B,A)=\sup_{x\in B}\di(x,A)$.\ed

\bp\label{p3.5} Let $\cD$ be an inclusion-closed collection of random subsets of $X$ and $\Phi$ a continuous RDS on $X$ over $(\W,\cF,\bP,(\theta_{t})_{t\in\R})$.
Suppose that $\{K(\w)\}_{\w\in\W}$ is a closed random absorbing set for $\Phi$ in $\cD$ and $\Phi$ is pullback $\cD$-asymptotically compact in $X$.
Then $\Phi$ has a unique $\cD$-random attractor $\{\sA(\w)\}_{\w\in\W}$ given by
$$\sA(\w)=\Cap_{T\geqslant0}\ol{\Cup_{t\geqslant T}\Phi(t,\theta_{-t}\w, K(\theta_{-t}\w))}.$$
\ep

\bp\label{p3.6} Let $\cD$ be an inclusion-closed collection of random subsets of $X$.
Given $\epsilon>0$, suppose that $\Phi_\epsilon$ is a random dynamical system over a measurable system $(\W,\cF,\bP,(\theta_t)_{t\in\R})$
which has a $\cD$-random attractor $\{\sA(\w)\}_{\w\in\W}$ and that $\Phi_0$ is a deterministic dynamical system defined on $X$ possessing a global attractor $\sA_0$.
Assume that the following conditions are satisfied,
\benu\item[(i)] for $\bP$-a.s. $\w\in\W$, $t\geqslant0$, $\epsilon_n\ra0$ and $x_n,\,x\in X$ with $x_n\ra x$, there holds that
$$\lim_{n\ra\8}\Phi_{\epsilon_n}(t,\w,x_n)=\Phi_0(t,x),$$
\item[(ii)] every $\Phi_\epsilon$ has a random absorbing set $\{K_\epsilon(\w)\}_{\w\in\W}\in\cD$ such that for some deterministic positive constant $C$ and for $\bP$-a.s. $\w\in\W$,
$$\limsup_{\epsilon\ra0}\|K_\epsilon(\w)\|_X\leqslant C,\hs\mb{and}$$
\item[(iii)] there exists $\epsilon_0>0$ such that for $\bP$-a.s. $\w\in\W$, the random set
$$\Cup_{0<\epsilon\leqslant\epsilon_0}\sA_\epsilon(\w)\hs \mb{is precompact in }X.$$
\eenu
Then for $\bP$-a.s. $\w\in\W$,
$${\rm dist}_X(\sA_\epsilon(\w),\sA_0)\ra0,\hs\mb{as }\epsilon\ra0.$$\ep

\subsection{Estimates for the solutions}\label{ss3.2}

In this subsection, we mainly present some uniform estimates on the solutions,
which in turn guarantees the existence of the random attractor.

In accord with the discussion in the last section, we know that given every initial datum $u_0\in H$ and $\w\in\W$,
the problem \eqref{1.1} -- \eqref{1.3} can generate a random dynamical system $\Phi_\epsilon$ via $u(t;\w,u_0)$ for each $\epsilon>0$,
where $\Phi_{\epsilon}:\R^+\X\W\X H\ra H$ is defined by
\be\label{3.1}
\Phi_\epsilon(t,\w,u_0)=u(t;\w,u_0),\hs\mb{for every }(t,\w,u_0)\in\R^+\X\W\X H.
\ee
Then $\Phi_\epsilon$ is a continuous random dynamical system over $(\W,\cF,\bP,(\theta_t)_{t\in\R})$ in $H$.

Suppose that $\widehat B=\{B(\w)\}_{\w\in\W}$ is a tempered family of bounded nonempty subsets of $H$, that is, for every $c>0$, $\w\in\W$,
\be\label{3.2}
\lim_{s\ra\8}\me^{-cs}\|B(\theta_{-s}\w))\|=0.
\ee
In the following, we always let $\cD$ be the collection of all tempered sets in $H$, that is, $\widehat B\in\cD$ if and only if the equation \eqref{3.2} holds.

We define below some random variables, which will be used frequently in the sequel.
Let
\be\label{3.3}m_{1\epsilon}(\w)=\me^{2\epsilon\max_{s\in[-1,0]}|z(\theta_s\w)|},\hs m_\epsilon(\w)=\me^{2\max_{s\in[-2,0]}|z(\theta_s\w)|}\ee
\be\label{3.4}\mb{and}\hs M_\epsilon(\w)=\int^0_{-\8}\me^{2s+2\epsilon\left|\int_{s}^{0} z(\theta_{\vsig}\w)\di\vsig\right|+2\epsilon|z(\theta_s\w)|}\di s.\ee
We see by \eqref{2.5} that $m_{1\epsilon}(\w)$, $m_\epsilon(\w)$ and $M_\epsilon(\w)$ are all tempered, and moreover for $\bP$-a.s. $\w\in\W$,
\be\label{3.5}1\leqslant m_{1\epsilon}(\w)\leqslant m_\epsilon(\w)\hs\mb{ and }\hs M_\epsilon(\w)\leqslant M_1(\w)\mb{ for all }\epsilon\in(0,1].\ee

Now we are ready to study the two specific cases for the kernel $G$.
In the following estimation, we allow $\epsilon$ to be zero and always let $\widehat B=\{B(\w)\}_{\w\in\W}$.
In the light of the transformational relations \eqref{2.6} and \eqref{2.13} between the solutions $u$ and $v$,
and the temperedness of $m_{1\epsilon}$, $m_{\epsilon}$ and $M_\epsilon$,
it is sufficient for us to make estimates of the solution $v(t;\w,v_0(\w))$.

\subsubsection{Estimates on solutions for the case of positive kernel}

Firstly, we give uniform estimates for the solutions to \eqref{2.7} -- \eqref{2.9} with $G$ satisfying ({\bf Gp}),
for which, we mainly refer to the similar estimating methods in \cite{G15}.

With the term $h$ involved, many details of the related estimations are not so clear in comparison with those in \cite{G15},
and these estimates are also pivotal for the argument for the upper semi-continuity of random attractors.
Therefore, we present a thorough estimation for our situation.

\bl\label{le3.7}Let $0\leqslant\epsilon\leqslant1$, $h\in H$ and $G$ satisfy ({\bf Gp}).
Assume that $v_0\in\widehat B$ with $\widehat B\in\cD$.
Then there exist a $t_{\widehat B}(\w)>0$ and a random variable $\rho_1(\w)>0$ such that, for all $\bP$-a.s. $\w\in\W$ and $t>t_{\widehat B}(\w)$, we have
\be\label{3.6}\|v(t;\theta_{-t}\w,v_0(\theta_{-t}\w))\|\lesssim\rho_1(\epsilon,\w),\ee
where
\be\label{3.7}\rho_1^2(\epsilon,\w)=1+M_\epsilon(\w)\(1+\|h\|^2\).\ee
\el
\bo\hs From \eqref{2.7}, the following equation formally holds for each $t\geqslant 0$ and $\w\in\W$,
\be\label{3.8}
\frac{1}{2}\frac{\di}{\di t}\|v\|^2+(a-\epsilon z(\theta_t\w))\|v\|^2+\|\De v\|^2+2(\De v,v)+\me^{2\epsilon z(\theta_t\w)}(vG*v^2,v)+\me^{-\epsilon z(\theta_t\w)}(h,v)=0.
\ee
By Schwarz's inequality and \eqref{2.10}, we have
\be\label{3.9} 2|(\De v,v)|\leqslant\frac12\|\De v\|^2+2\|v\|^2,\hs(vG*v^2,v)\geqslant\al\|v\|^4\ee
\be\label{3.10}\mb{and}\hs\me^{-\epsilon z(\theta_t\w)}|(h,v)|\leqslant\frac{1}{2}\me^{-2\epsilon z(\theta_t\w)}\|h\|^2+\frac12\|v\|^2.\ee
Combining \eqref{3.8} -- \eqref{3.10}, we have
\begin{align*}&\frac{\di}{\di t}\|v\|^2+2(1-\epsilon z(\theta_t\w))\|v\|^2+\|\De v\|^2+2\al\me^{2\epsilon z(\theta_t\w)}\|v\|^4\\
\leqslant&(7-2a)\|v\|^2+\me^{-2\epsilon z(\theta_t\w)}\|h\|^2\\
\leqslant&2\al\me^{2\epsilon z(\theta_t\w)}\|v\|^4+c\me^{-2\epsilon z(\theta_t\w)}(1+\|h\|^2),\end{align*}
$$\mb{and}\hs\frac{\di}{\di t}\|v\|^2+2(1-\epsilon z(\theta_t\w))\|v\|^2+\|\De v\|^2\leqslant c\me^{-2\epsilon z(\theta_t\w)}\(1+\|h\|^2\).$$
By Gronwall's lemma, we have for $T\geqslant\sig\geqslant0$,
\begin{align}
&\|v(T;\w,v_0(\w))\|^2+\int_{\sig}^{T}\me^{-2(T-s)+2\epsilon\int_s^{T} z(\theta_{\vsig}\w)\di\vsig}\|\De v\|^2\di s\notag\\
\leqslant&\me^{-2(T-\sig)+2\epsilon\int_{\sig}^{T} z(\theta_{\vsig}\w)\di\vsig}\|v(\sig;\w,v_0(\w))\|^2\label{3.11}\\
&+c\(1+\|h\|^2\)\int_{\sig}^{T}\me^{-2(T-s)+2\epsilon\int_s^{T} z(\theta_{\vsig}\w)\di\vsig-2\epsilon z(\theta_{s}\w)}\di s.\notag\end{align}
Replacing $\sig$ and $T$ in \eqref{3.11} by $0$ and $t$ respectively, we obtain
\begin{align}
&\|v(t;\w,v_0(\w))\|^2+\int_{0}^{t}\me^{-2(t-s)+2\epsilon\int_s^t z(\theta_{\vsig}\w)\di\vsig}\|\De v\|^2\di s\notag\\
\leqslant&\me^{-2t+2\epsilon\int_0^t z(\theta_{\vsig}\w)\di\vsig}\|v_0(\w)\|^2\label{3.12}\\
&+c\(1+\|h\|^2\)\int_0^t\me^{-2(t-s)+2\epsilon\int_s^t z(\theta_{\vsig}\w)\di\vsig-2\epsilon z(\theta_{s}\w)}\di s.\notag\end{align}
By replacing $\w$ by $\theta_{-t}\w$ in \eqref{3.12}, we have
\begin{align}
&\|v(t;\theta_{-t}\w,v_0(\theta_{-t}\w))\|^2+\int_{-t}^0\me^{2s+2\epsilon\int_s^0 z(\theta_{\vsig}\w)\di\vsig}\|\De v(s+t;\theta_{-t}\w,v_0(\theta_{-t}\w))\|^2\di s\notag\\
\leqslant&\me^{-2t+2\epsilon\int_{-t}^{0} z(\theta_{\vsig}\w)\di\vsig}\|v_0(\theta_{-t}\w))\|^2+c\(1+\|h\|^2\)\int_{-t}^0\me^{2s+2\epsilon\int_{s}^0z(\theta_{\vsig}\w)\di\vsig-2\epsilon z(\theta_{s}\w)}\di s\label{3.13}\\
:=&E_\epsilon(t,\w).\notag\end{align}
Due to \eqref{2.5}, we know that for each $\gamma>0$, there is a positive random variable $R(\gamma,\w)$ such that for all $\bP$-a.s. $\w\in\W$ and $t\in\R$,
\be\label{3.14}\max\left\{|z(\theta_t\w)|,\,\left|\int^t_0z(\theta_\vsig\w)\di\vsig\right|\right\}\leqslant\gamma|t|+R(\gamma,\w).\ee
Then when $0\leqslant\epsilon\leqslant1$, by the fact that $v_0\in\widehat B$ and \eqref{3.2}, there is a positive number $t_{\widehat B}(\w)$ independent of $\epsilon$, such that as $t>t_{\widehat B}(\w)$,
\be\label{3.15}\me^{-2t+2\epsilon\int^{0}_{-t}z(\theta_{\vsig}\w)\di\vsig}\|v_0(\theta_{-t}\w)\|^2\leqslant1,\ee
Therefore by \eqref{3.4}, we obtain that when $t>t_{\widehat B}(\w)$,
\be\label{3.16}E_\epsilon(t,\w)\lesssim\rho_1^2(\epsilon,\w),\ee
where $\rho_1^2(\epsilon,\w)$ is defined as \eqref{3.7}.
Thus \eqref{3.6} follows from \eqref{3.13} and \eqref{3.16}.\eo

\bl\label{le3.8} Let $0\leqslant\epsilon\leqslant 1$, $h\in H$ and $G$ satisfy ({\bf Gp}).
Assume that $v_0\in\widehat B$ with $\widehat B\in\cD$.
Then when $t>t_{\widehat B}(\w)$, for all $\bP$-a.s. $\w\in\W$,
\be\label{3.17}\|\De v(t;\theta_{-t}\w,v_0(\theta_{-t}\w))\|\lesssim\rho_2(\epsilon,\w),\ee
where
\be\label{3.18}\rho_2^2(\epsilon,\w)=m_{1\epsilon}^2(\w)\left[1+M_\epsilon(\w)\(1+\|h\|^2\)\right]^2.\ee
\el

\bo Multiplying \eqref{2.7} by $\De^2 v$ and integrating it over $U$, we have for $t\geqslant 0$ and $\w\in\W$,
\be\label{3.19}
\ba{rl}&\disp\frac{1}{2}\frac{\di}{\di t}\|\De v\|^2+(a-\epsilon z(\theta_t\w))\|\De v\|^2+\|\De^2 v\|^2\\
=&-2(\De v,\De^2v)-\me^{-\epsilon z(\theta_t\w)}(h,\De^2 v)-\me^{2\epsilon z(\theta_t\w)}(vG*v^2,\De^2v).\ea
\ee
By Schwartz's inequality and embedding inequalities, we have
\be 2|(\De v,\De^2v)|\leqslant4\|\De v\|^2+\frac14\|\De^2 v\|^2,\ee
\be\me^{-\epsilon z(\theta_t\w)}|(h,\De^2v)|\leqslant\me^{-2\epsilon z(\theta_t\w)}\|h\|^2+\frac14\|\De^2v\|^2\hs\mb{and}\ee
\begin{align}
&|(vG*v^2,\De^2v)|=|(\De vG*v^2+2\na v\na G*v^2+v\De G*v^2,\De v)|\notag\\
\leqslant&\b\|v\|^2\int_U(|\De v|+2|\na v|+|v|)|\De v|\di x\leqslant c\|v\|^2\|\De v\|^2.\label{3.22}
\end{align}
Combining \eqref{3.19} -- \eqref{3.22}, we obtain the following differential inequality,
\be\label{3.23}\ba{rl}&\disp\frac{\di}{\di t}\|\De v\|^2+2(1-\epsilon z(\theta_t\w))\|\De v\|^2+\|\De^2v\|^2\\
\leqslant&c\|\De v\|^2+c\me^{2\epsilon z(\theta_t\w)}\|v\|^2\|\De v\|^2+c\me^{-2\epsilon z(\theta_t\w)}\|h\|^2.\ea\ee
Applying Gronwall's lemma to \eqref{3.23}, we have that
\begin{align}
&\|\De v(T;\w,v_0(\w))\|^2+\int_\sig^T\me^{-2(T-s)+2\epsilon\int_s^Tz(\theta_\vsig\w)\di\vsig}\|\De^2 v\|^2\di s\notag\\
\leqslant&\me^{-2(T-\sig)+2\epsilon\int_\sig^Tz(\theta_\vsig\w)\di\vsig}\|\De v(\sig;\w,v_0(\w))\|^2\label{3.24}\\
&+c\int_\sig^T\me^{-2(T-s)+2\epsilon\int_s^Tz(\theta_\vsig\w)\di\vsig}\(\me^{-2\epsilon z(\theta_s\w)}\|h\|^2+\|\De v\|^2+\me^{2\epsilon z(\theta_s\w)}\|v\|^2\|\De v\|^2\)\di s.\notag\end{align}
Replacing $T$ and $\w$ by $t+1$ and $\theta_{-t-1}\w$ respectively and integrating \eqref{3.24} from $t$ to $t+1$ with respect to $\sig$, and by \eqref{3.13}, we have that
\begin{align}&
  \|\De v(t+1;\theta_{-t-1}\w,v_0(\theta_{-t-1}\w))\|^2\notag\\
\leqslant&
  \int_{-1}^0\me^{2s+2\epsilon\int_s^0z(\theta_{\vsig}\w)\di\vsig}\|\De v(s+t+1;\theta_{-t-1}\w,v_0(\theta_{-t-1}\w))\|^2\di s\notag\\
&
  +c\|h\|^2\int_{-1}^0\int_{\sig}^{0}\me^{2s+2\epsilon\int_s^0z(\theta_\vsig\w)\di\vsig-2\epsilon z(\theta_s\w)}\di s\di\sig\notag\\
&
  +c\int_{-1}^0\int_{\sig}^{0}\me^{2s+2\epsilon\int_s^0z(\theta_\vsig\w)\di\vsig}\|\De v(s+t+1;\theta_{-t-1}\w,v_0(\theta_{-t-1}\w))\|^2\di s\di\sig\notag\\
&
  +c\int_{-1}^0\int_{\sig}^{0}\me^{2s+2\epsilon\int_s^0z(\theta_\vsig\w)\di\vsig+2\epsilon z(\theta_s\w)}\|v\|^2\|\De v\|^2\di s\di\sig\notag\\
\lesssim&
  M_\epsilon(\w)\|h\|^2+E_\epsilon(t+1,\w)+\int_{-1}^0\me^{2s+2\epsilon\int_s^0z(\theta_\vsig\w)\di\vsig+2\epsilon z(\theta_s\w)}\|v\|^2\|\De v\|^2\di s.\label{3.25}
\end{align}
Replacing $t$ and $\w$ by $s+t+1$ and $\theta_{-t-1}\w$ respectively in \eqref{3.12} with $s\in[-t-1,0]$ and by \eqref{3.13}, we have
\begin{align}&
  \me^{2s+2\epsilon\int_s^0z(\theta_\vsig\w)\di\vsig}\|v(s+t+1;\theta_{-t-1}\w,v_0(\theta_{-t-1}\w))\|^2\notag\\
\leqslant&
  \me^{-2(t+1)+2\epsilon\int_{-t-1}^0z(\theta_\vsig\w)\di\vsig}\|v_0(\theta_{-t-1}\w)\|^2+c\(1+\|h\|^2\)\int_{-t-1}^0\me^{2r
  +2\epsilon\int_{r}^0z(\theta_\vsig\w)\di\vsig-2\epsilon z(\theta_{r}\w)}\di r\notag\\
\leqslant&
  E_\epsilon(t+1,\w).\label{3.26}
\end{align}
Hence
\begin{align}&
  \int_{-1}^0\me^{2s+2\epsilon\int_s^0z(\theta_\vsig\w)\di\vsig+2\epsilon z(\theta_s\w)}\|v\|^2\|\De v\|^2\di s\notag\\
\leqslant&
  E_\epsilon(t+1,\w)\int_{-1}^0\me^{-2s-2\epsilon\int_s^0z(\theta_\vsig\w)\di\vsig+2\epsilon z(\theta_s\w)}\me^{2s
  +2\epsilon\int_s^0z(\theta_\vsig\w)\di\vsig}\|\De v\|^2\di s\notag\\
\leqslant&
  \me^2m_{1\epsilon}^2(\w)E^2_\epsilon(t+1,\w).
  \label{3.27}
\end{align}
Combining \eqref{3.25}, \eqref{3.27}, \eqref{3.16}, \eqref{3.7} and \eqref{3.5}, we obtain \eqref{3.17} when $t>t_{\widehat B}(\w)$.
The proof is now complete.
\eo

Since $m_{1\epsilon}(\w)$ and $M_\epsilon(\w)$ are both tempered, so are the random variables $\rho_i(\w)$ for $i=1,2$.
Lemma \ref{le3.7} guarantees the existence of the random absorbing set in $\cD$ for the system generated by \eqref{2.7} -- \eqref{2.9}.
Lemma \ref{le3.8} gives the pullback $\cD$-asymptotic compactness in $H$ due to the compact embedding from $V$ into $H$.
Then combining these results and \eqref{2.11}, we obtain the existence and regularity of the random attractor for \eqref{2.7} -- \eqref{2.9} (see \cite{G15}),
which is stated as follows in details.

\bt\label{th3.9}Let $0<\epsilon\leqslant1$.
There is a $\cD$-random attractor contained in $V$ for the random system generated by \eqref{2.7} -- \eqref{2.9}.
\et

Theorem \ref{th3.9} motivates us to study the upper-continuity of random attractors in $V$.
To this goal, we still need to show the pullback $\cD$-asymptotic compactness of $\Phi_\epsilon$ in $V$.
And it is necessary to work in the framework of fractional spaces $\sD(A^{\mu})$ for $\mu\in(\frac12,1)$.

By the settings in Section \ref{s2}, the problem \eqref{2.7} -- \eqref{2.9} can be rewritten in an abstract form,
\be\label{3.28}\left\{\ba{l}\disp\frac{\di v}{\di t}=-(A+1-\epsilon z(\theta_t\w))v+f_\epsilon(\theta_t\w,v),\; t>0,\\
v(0)=v_0(\w)=\me^{-\epsilon z(\w)}u_0\in H,\hs\mb{where}\\
f_\epsilon(\w,v)=2A^{\frac12}v-(a-1)v-\me^{2\epsilon z(\w)}vG*v^2-\me^{-\epsilon z(\w)}h.\ea\right.\ee
The solution to \eqref{3.28} satisfies the following integral equation (see \cite{Hen81}), for all $t\geqslant s\geqslant0$,
\begin{align}v(t;\w,v_0(\w))=&\me^{-A(t-s)}\me^{-(t-s)+\epsilon\int_s^tz(\theta_\vsig\w)\di\vsig}v(s;\w,v_0(\w))\notag\\
&+\int_s^t\me^{-A(t-\sig)}\me^{-(t-\sig)+\epsilon\int_\sig^tz(\theta_\vsig\w)\di\vsig}f_\epsilon(\theta_\sig\w,v(\sig;\w,v_0(\w)))\di\sig.\label{3.29}\end{align}
We have the following estimate for the $\sD(A^\mu)$-norm of the solution $v(t;\w,v_0(\w))$.

\bl\label{le3.9} Let $0\leqslant\epsilon\leqslant1$, $h\in H$ and $G$ satisfy ({\bf Gp}).
Assume that $v_0\in\widehat B$ with $\widehat B\in\cD$ and $\mu\in(\frac12,1)$.
Then when $t>t_{\widehat B}(\w)$, for all $\bP$-a.s. $\w\in\W$,
\be\label{3.30}\|v(t;\theta_{-t}\w,v_0(\theta_{-t}\w))\|_\mu\lesssim\rho_3(\epsilon,\w),\ee
where
\be\label{3.31}\rho_3(\epsilon,\w):=m^5_\epsilon(\w)\left[1+M_\epsilon(\w)\(1+\|h\|^2\)\right]^2.\ee
\el
\bo By replacement of $\sig$, $T$ and $\w$ by $t$, $t+1$ and $\theta_{-t-1}\w$, respectively, in \eqref{3.24}, we have
\begin{align*}&
  F_\epsilon(t+1,\w)\\
:=&
  \int_{-1}^{0}\me^{2s+2\epsilon\int_s^0z(\theta_\vsig\w)\di\vsig}\|\De^2 v(s+t+1;\theta_{-t-1}\w,v_0(\theta_{-t-1}\w))\|^2\di s\\
\leqslant&
  \me^{-2+2\epsilon\int^0_{-1}z(\theta_\vsig\w)\di\vsig}\|\De v(t;\theta_{-t-1}\w,v_0(\theta_{-t-1}\w))\|^2\\
&
  +c\int_{-1}^0\me^{2s+2\epsilon\int_s^0z(\theta_\vsig\w)\di\vsig}\(\me^{-2\epsilon z(\theta_s\w)}\|h\|^2+\|\De v\|^2+\me^{2\epsilon z(\theta_s\w)}\|v\|^2\|\De v\|^2\)\di s.
\end{align*}
Then similar to the discussions of \eqref{3.17} and \eqref{3.25}, since $v_0\in\widehat B$, then when $t>t_{\widehat B}(\w)$,
\be\label{3.32}F_\epsilon(t+1,\w)\leqslant m_{1\epsilon}(\w)\rho_2^2(\epsilon,\theta_{-1}\w)+\rho_2^2(\epsilon,\w).\ee
By the definitions \eqref{3.3} and \eqref{3.4}, we have $m_{1\epsilon}(\theta_{-1}\w)\leqslant m_\epsilon(\w)$ and
\begin{align}
  M_\epsilon(\theta_{-1}\w)
=&
  \int^0_{-\8}\me^{2s+2\epsilon\left|\int_{s}^{0} z(\theta_{\vsig}\theta_{-1}\w)\di\vsig\right|+2\epsilon|z(\theta_s\theta_{-1}\w)|}\di s\notag\\
=&
  \int^{-1}_{-\8}\me^{2(s+1)+2\epsilon\left|\(\int_{s}^0-\int_{-1}^0\)z(\theta_{\vsig}\w)\di\vsig\right|+2\epsilon|z(\theta_s\w)|}\di s\notag\\
\leqslant&
  \me^2m_{1\epsilon}(\w)M_\epsilon(\w).\label{3.33}
\end{align}
It follows from \eqref{3.18}, \eqref{3.32}, \eqref{3.33} and \eqref{3.5} that when $t>t_{\widehat B}(\w)$,
\be\label{3.34}F_\epsilon(t+1,\w)\lesssim m_{\epsilon}^5(\w)\left[1+M_\epsilon(\w)(1+\|h\|^2)\right]^2.\ee

Fix $\mu\in(\frac12,1)$.
Then for $t>0$, by \eqref{2.11}, \eqref{2.1}, \eqref{2.2} and the embeddings of spaces of fractional order (\cite{Hen81}), we have that for $0<s\leqslant t$,
\be\label{3.35}\|A^\mu\me^{-A(t-s)}v(s)\|=\|\me^{-A(t-s)}A^\mu v(s)\|\lesssim\me^{-\lam_0(t-s)}\|\De^2v(s)\|,\ee
\be\|A^\mu\me^{-A(t-s)}A^{\frac12}v(s)\|=\|A^{\mu-\frac12}\me^{-A(t-s)}Av(s)\|\lesssim(t-s)^{\frac12-\mu}\me^{-\lam_0(t-s)}\|\De^2v(s)\|,\ee
\be\|A^{\mu}\me^{-A(t-s)}v(s)\|\lesssim\me^{-\lam_0(t-s)}\|\De^2v(s)\|,\ee
\be\|A^{\mu}\me^{-A(t-s)}h\|\lesssim(t-s)^{-\mu}\me^{-\lam_0(t-s)}\|h\|\ee
and by the calculation in \eqref{3.22} additionally,
\begin{align}&
  \|A^{\mu}\me^{-A(t-s)}\(v(s)G*v^2(s)\)\|\notag\\
=&
  \|A^{\mu-\frac12}\me^{-A(t-s)}\De\(v(s)G*v^2(s)\)\|\notag\\
\lesssim&
  (t-s)^{\frac12-\mu}\me^{-\lam_0(t-s)}\|\De v(s)\|\|v(s)\|^2,\label{3.39}
\end{align}
where we have used the embeddings of Sobolev spaces.

Consider the $\sD(A^\mu)$-norm of \eqref{3.29} and use the inequalities \eqref{3.35} -- \eqref{3.39}.
It can be seen that
\begin{align}&
  \|v(t;\w,v_0(\w))\|_\mu\notag\\
\lesssim&
  \me^{-(t-s)+\epsilon\int_s^tz(\theta_\vsig\w)\di\vsig}\|\De^2v(s;\w,v_0(\w))\|\notag\\
&
  +\int_s^t\me^{-(t-\sig)+\epsilon\int_\sig^tz(\theta_\vsig\w)\di\vsig}\left[(t-\sig)^{\frac12-\mu}+|a-1|\right]\|\De^2v(\sig)\|\di\sig\notag\\
&
  +\int_s^t\me^{-(t-\sig)+\epsilon\int_\sig^tz(\theta_\vsig\w)\di\vsig-\epsilon z(\theta_\sig\w)}(t-\sig)^{-\mu}\|h\|\di\sig\notag\\
&
  +\int_s^t\me^{-(t-\sig)+\epsilon\int_\sig^tz(\theta_\vsig\w)\di\vsig+2\epsilon z(\theta_\sig\w)}(t-\sig)^{\frac12-\mu}\|\De v(\sig)\|\|v(\sig)\|^2\di\sig.\label{3.40}
\end{align}
Replacing $t$ and $\w$ by $t+1$ and $\theta_{-t-1}\w$ respectively in \eqref{3.40}, and integrating it over $[t,t+1]$ with respect to $s$, we obtain
\begin{align}&\|v(t+1;\theta_{-t-1}\w,v_0(\theta_{-t-1}\w))\|_\mu\notag\\
\lesssim&
  \int_t^{t+1}\me^{-(t+1-s)+\epsilon\int_s^{t+1}z(\theta_{\vsig-t-1}\w)\di\vsig}\|\De^2v(s;\theta_{-t-1}\w,v_0(\theta_{-t-1}\w))\|\di s\notag\\
&
  +\int_t^{t+1}\me^{-(t+1-s)+\epsilon\int_s^{t+1}z(\theta_{\vsig-t-1}\w)\di\vsig}\left[(t+1-s)^{\frac12-\mu}+|a-1|\right]\|\De^2v(s)\|\di s\notag\\
&
  +\int_t^{t+1}\me^{-(t+1-s)+\epsilon\int_s^{t+1}z(\theta_{\vsig-t-1}\w)\di\vsig-\epsilon z(\theta_{s-t-1}\w)}(t+1-s)^{-\mu}\|h\|\di s\notag\\
&
  +\int_t^{t+1}\me^{-(t+1-s)+\epsilon\int_s^{t+1}z(\theta_{\vsig-t-1}\w)\di\vsig+2\epsilon z(\theta_{s-t-1}\w)}(t+1-s)^{\frac12-\mu}\|\De v(s)\|\|v(s)\|^2\di s\notag\\
\lesssim&
  \int_{-1}^0\me^{s+\epsilon\int_s^0z(\theta_\vsig\w)\di\vsig}\left[(-s)^{\frac12-\mu}+1\right]\|\De^2v(s+t+1;\theta_{-t-1}\w,v_0(\theta_{-t-1}\w))\|\di s\notag\\
&
 +\|h\|\int_{-1}^0\me^{s+\epsilon\int_s^0z(\theta_\vsig\w)\di\vsig-\epsilon z(\theta_s\w)}(-s)^{-\mu}\di s\notag\\
&
 +\int_{-1}^0\me^{s+\epsilon\int_s^0z(\theta_\vsig\w)\di\vsig+2\epsilon z(\theta_s\w)}(-s)^{\frac12-\mu}\|\De v(s+t+1)\|\|v(s+t+1)\|^2\di s.\label{3.41}
\end{align}
When $t>t_{\widehat B}(\w)$, applying \eqref{3.26}, \eqref{3.27} and Schwarz's inequality to \eqref{3.41}, we have
\begin{align*}&
  \|v(t+1,\theta_{-t-1}\w,v_0(\theta_{-t-1}\w))\|_\mu\\
\lesssim&
  F_\epsilon(t+1,\w)+\int_{-1}^0\left[(-s)^{\frac12-\mu}+1\right]^2\di s+m_{1\epsilon}\|h\|\int_{-1}^0\me^s(-s)^{-\mu}\di s\\
&
  +\int_{-1}^0\me^{2s+2\epsilon\int_s^0z(\theta_{\vsig}\w)\di\vsig+4\epsilon z(\theta_{s}\w)}\|v(s+t+1)\|^2\|\De v(s+t+1)\|^2\di s\\
&
  +\int_{-1}^0(-s)^{1-2\mu}\|v(s+t+1)\|^2\di s\\
\lesssim&
  F_\epsilon(t+1,\w)+\me^2m^3_{1\epsilon}E_\epsilon^2(t+1,\w)+\int_{-1}^0\left[(-s)^{\frac12-\mu}+1\right]^2\di s\\
&
  +m_{1\epsilon}(\w)\|h\|\int_{-1}^0\me^s(-s)^{-\mu}\di s+\me^2m_{1\epsilon}(\w)E_\epsilon(t+1,\w)\int_{-1}^0(-s)^{1-2\mu}\di s.
\end{align*}
Hence by \eqref{3.5}, \eqref{3.16} and \eqref{3.34}, we see that when $t>t_{\widehat B}(\w)$, the conclusion \eqref{3.30} holds.
\eo
\br The pullback $\cD$-asymptotic compactness of $\Phi_\epsilon$ in $V$ follows immediately from Lemma \ref{le3.9}, \eqref{2.11}
and the compact embedding from $\sD(A^\mu)$ into $V$ for $\mu\in(\frac12,1)$.
Indeed, Lemmas \ref{le3.8} and \ref{le3.9} indicate that the random attractor given by Theorem \ref{th3.9} is also a random attractor
in the sense of $V$-norm.
\er

\subsubsection{Estimates on solutions for the case of non-negative kernel}

For the non-negative kernels, the analysis above does not work smoothly here.
However, if we adopt the special non-negative kernel $G=J_{\de_0}$, or saying, under the assumption ({\bf Gn}),
we can also obtain similar estimates on the solutions to \eqref{2.7} -- \eqref{2.9} as follows by a slight modification in the proofs.

\bl\label{le3.11}Let $0\leqslant\epsilon\leqslant1$, $h\in H$ and $G$ satisfy ({\bf Gn}).
Assume that $v_0\in\widehat B$ with $\widehat B\in\cD$.
Then there exists a $t'_{\widehat B}(\w)>0$ such that, for all $\bP$-a.s. $\w\in\W$ and $t>t'_{\widehat B}(\w)$, we have
\be\label{3.42}\|v(t;\theta_{-t}\w,v_0(\theta_{-t}\w))\|\lesssim\rho_1(\epsilon,\w),\ee
\be\label{3.43}\|\De v(t;\theta_{-t}\w,v_0(\theta_{-t}\w))\|\lesssim\rho_2(\epsilon,\w)\ee
\be\label{3.44}\mb{and}\hs\|v(t;\theta_{-t}\w,v_0(\theta_{-t}\w))\|_\mu\lesssim\rho_3(\epsilon,\w),\ee
where $\mu\in(\frac12,1)$, $\rho_i$ is defined in \eqref{3.7}, \eqref{3.18} and \eqref{3.31}, respectively.
In particular, when $0<\epsilon\leqslant1$, there is a $\cD$-random attractor contained in $V$ for the random system generated by \eqref{2.7} -- \eqref{2.9}.
\el
\br It is necessary to mention that the $\lesssim$'s given in \eqref{3.42} -- \eqref{3.44} imply different positive constants from those in \eqref{3.6}, \eqref{3.17} and \eqref{3.30}, respectively.\er

\noindent{\it Proof of Lemma \ref{le3.11}.}
The key difference between the estimations in the two cases lies mainly in those for $\|v(t;\theta_{-t}\w,v_0(\theta_{-t}\w))\|$ presented in Lemma \ref{le3.7}.
We hence pay the most attention to this estimation.

From \eqref{2.7}, we similarly have that for each $t\geqslant 0$ and $\w\in\W$,
the equation \eqref{3.8} and the inequalities \eqref{3.9} and \eqref{3.10} hold with only the estimate for $(vG*v^2,v)$ replaced by
\be\label{3.45}(vG*v^2,v)\geqslant\al c_0\|v\|^4-\de\|v\|^2,\ee
where $c_0$ is obtained from the embedding from $L^4(U)$ into $H$.
Then similarly we have
\begin{align*}&\frac{\di}{\di t}\|v\|^2+2(1-\epsilon z(\theta_t\w))\|v\|^2+\|\De v\|^2+2\al c_0\me^{2\epsilon z(\theta_t\w)}\|v\|^4\\
\leqslant&(7-2a)\|v\|^2+\me^{-2\epsilon z(\theta_t\w)}\|h\|^2+2\de\me^{2\epsilon z(\theta_t\w)}\|v\|^2\\
\leqslant&2\al c_0\me^{2\epsilon z(\theta_t\w)}\|v\|^4+c\me^{-2\epsilon z(\theta_t\w)}(1+\|h\|^2)+c\me^{2\epsilon z(\theta_t\w)}\end{align*}
$$\mb{and}\hs\frac{\di}{\di t}\|v\|^2+2(1-\epsilon z(\theta_t\w))\|v\|^2+\|\De v\|^2\leqslant c\me^{-2\epsilon z(\theta_t\w)}\(1+\|h\|^2\)+c\me^{2\epsilon z(\theta_t\w)},$$
By Gronwall's lemma, we have for $T\geqslant\sig$,
\begin{align}
&\|v(T;\w,v_0(\w))\|^2+\int_{\sig}^{T}\me^{-2(T-s)+2\epsilon\int_s^{T} z(\theta_{\vsig}\w)\di\vsig}\|\De v\|^2\di s\notag\\
\leqslant&\me^{-2(T-\sig)+2\epsilon\int_{\sig}^{T} z(\theta_{\vsig}\w)\di\vsig}\|v(\sig;\w,v_0(\w))\|^2\label{3.46}\\
&+c\(1+\|h\|^2\)\int_{\sig}^{T}\me^{-2(T-s)+2\epsilon\int_s^{T} z(\theta_{\vsig}\w)\di\vsig-2\epsilon z(\theta_{s}\w)}\di s\notag\\
&+c\int_{\sig}^{T}\me^{-2(T-s)+2\epsilon\int_s^{T} z(\theta_{\vsig}\w)\di\vsig+2\epsilon z(\theta_{s}\w)}\di s.\notag\end{align}
Replacing $\sig$, $T$ and $w$ in \eqref{3.46} by $0$, $t$ and $\theta_{-t}\w$ respectively, we obtain
\begin{align*}&
  \|v(t;\theta_{-t}\w,v_0(\theta_{-t}\w))\|^2+\int_{-t}^0\me^{2s+2\epsilon\int_s^0 z(\theta_{\vsig}\w)\di\vsig}\|\De v(s+t;\theta_{-t}\w,v_0(\theta_{-t}\w))\|^2\di s\\
\leqslant&
  \me^{-2t+2\epsilon\int_{-t}^{0} z(\theta_{\vsig}\w)\di\vsig}\|v_0(\theta_{-t}\w))\|^2+c\(1+\|h\|^2\)\int_{-t}^0\me^{2s+2\epsilon\int_{s}^0z(\theta_{\vsig}\w)\di\vsig-2\epsilon z(\theta_{s}\w)}\di s\\
&
  +c\int_{-t}^0\me^{2s+2\epsilon\int_{s}^0z(\theta_{\vsig}\w)\di\vsig+2\epsilon z(\theta_{s}\w)}\di s\\
:=&
  E'_\epsilon(t,\w).
\end{align*}
Then when $0\leqslant\epsilon\leqslant1$, by the fact that $v_0\in\widehat{B}$ and \eqref{3.14}, there is a positive number $t'_{\widehat{B}}(\w)$ independent of $\epsilon$, such that as $t>t'_{\widehat{B}}(\w)$, \eqref{3.15} is valid.

As a result, when $t>t'_{\widehat{B}}(\w)$,
$$E'_\epsilon(t,\w)\leqslant1+c\left[(1+\|h\|^2)M_\epsilon(\w)+M_\epsilon(\w)\right]\lesssim 1+M_\epsilon(1+\|h\|^2),$$
from which \eqref{3.42} follows.

The estimates in \eqref{3.43} and \eqref{3.44} follow from the same procedures of the proofs of Lemmas \ref{le3.8} and \ref{le3.9}.
We thus omit them.
\qed

\subsection{Upper semi-continuity of random attractors}\label{s4}

To indicate the dependence of solutions on $\epsilon$, we write the solution to the problem \eqref{1.1} -- \eqref{1.3} as $u_\epsilon$.
And we denote the solution and the semigroup of the following deterministic equation by $u$ and $\Phi_0$, respecively,
\be\label{3.47}
\frac{\pa u}{\pa t}+\De^2u+2\De u+au+uG*u^2+h(x)=0, \qquad x\in U,~ t>0.
\ee

\subsubsection{The case of positive kernel}

We first consider the nonlocal stochastic SHE with positive kernel $G$, i.e., $G$ satisfying \eqref{2.10}.
The existence of the global attractor of $\Phi_0$ in $V$ has been obtained by the discussion in Subsection \ref{ss3.2}.
We set $\widehat\sA_\epsilon=\{\sA_\epsilon(\w)\}_{\w\in\W}$ to be the $\cD$-random attractor of $\Phi_\epsilon$ and $\sA_0$ the global attractor of $\Phi_0$.

\bl\label{le3.13} Let $0<\epsilon\leqslant1$, $h\in H$ and $G$ satisfy ({\bf Gp}).
Then every $\Phi_\epsilon$ has a random absorbing set $\widehat K_\epsilon=\{K_\epsilon(\w)\}_{\w\in\W}\in\cD$ such that
for some deterministic positive constant $\ol C$ and all $\bP$-a.s. $\w\in\W$,
\be\label{3.48}\limsup_{\epsilon\ra0^+}\|K_\epsilon(\w)\|\leqslant\ol C.\ee
\el
\bo By Lemma \ref{le3.8}, \ref{le3.11} and \eqref{2.13}, for each $\epsilon\in(0,1]$, the process $\Phi_\epsilon$ has a closed and tempered random absorbing set $\widehat K_\epsilon=\{K_\epsilon(\w)\}_{\w\in\W}$ such that
$$K_\epsilon(\w)=\{u\in V:\|\De u(\w)\|\lesssim\me^{\epsilon|z(\w)|}\rho_2(\epsilon,\w)\},$$
in $V$, where $\rho_2$ is given in \eqref{3.18}.

Note that
\be\label{3.49}\|K_\epsilon(\w)\|^2\lesssim\me^{2\epsilon|z(\w)|}\rho^2_2(\epsilon,\w)\lesssim m_{1\epsilon}^3(\w)\left[1+M_\epsilon(\w)\(1+\|h\|^2\)\right]^2\ee
Fix an $\w\in\W$. By definition \eqref{3.3}, we easily have
\be\lim_{\epsilon\ra0^+}m_{1\epsilon}(\w)=1.\ee
Using \eqref{3.14} and taking $\gamma=1/4$ therein, we have
\begin{align}M_\epsilon(\w)\leqslant&\int_{-\8}^0\exp\(2s+4\epsilon\(\gam|s|+R\(\gam,\w\)\)\)\di s\notag\\
=&\me^{4\epsilon R\(\gam,\w\)}\int_{-\8}^0\me^{(2-\epsilon)s}\di s\notag\\
=&\frac1{2-\epsilon}\me^{4\epsilon R\(\gam,\w\)}\ra\frac12,\label{3.51}\end{align}
as $\epsilon\ra0^+$.
It follows from \eqref{3.49} -- \eqref{3.51} that
$$\limsup_{\epsilon\ra0^+}\|K_\epsilon(\w)\|\lesssim\(1+\frac{1+\|h\|^2}{2}\),$$
which indeed implies \eqref{3.48} and ends the proof.
\eo

\bl\label{le3.14} Let $0<\epsilon\leqslant1$, $h\in H$ and $G$ satisfy ({\bf Gp}).
Then for $\bP$-a.s. $\w\in\W$, the union $\cup_{0<\epsilon\leqslant1}\sA_\epsilon(\w)$ is precompact in $V$.\el
\bo Following the proof of Lemma \ref{le3.13}, it is easy to see that
$$\cup_{0<\epsilon\leqslant1}\sA_\epsilon(\w)\subset K_1(\w)\hs\mb{and}\hs \widehat K_1\in\cD.$$
Let $\mu\in(\frac12,1)$.
By Lemma \ref{le3.9}, there is a positive $t_{\widehat K_1}(\w)$ independent of $\epsilon$ such that when $t>t_{\widehat K_1}(\w)$, for $\bP$-a.s. $\w\in\W$,
$$\|\Phi_\epsilon(t,\theta_{-t}\w,K_1(\theta_{-t}\w))\|_{\mu}=\sup_{v_0\in \widehat K_1}\|\me^{\epsilon z(\w)}v(t;\theta_{-t}\w,v_0(\theta_{-t}\w))\|_\mu\lesssim\me^{\epsilon z(\w)}\rho_3(\epsilon,\w)\lesssim\me^{|z(\w)|}\rho_3(1,\w).$$
By the invariance of attractors, we have for $\bP$-a.s. $\w\in\W$,
$$\|\cup_{0<\epsilon\leqslant1}\sA_\epsilon(\w)\|_{\mu}\leqslant\|\cup_{0<\epsilon\leqslant1}\Phi_\epsilon(t,\theta_{-t}\w,K_1(\theta_{-t}\w))\|_{\mu}\lesssim\me^{|z(\w)|}\rho_3(1,\w).$$
By the compact embedding from $\sD(A^{\mu})$ into $V$, we know $\cup_{0<\epsilon\leqslant1}\sA_\epsilon(\w)$ is precompact in $V$, for $\bP$-a.s. $\w\in\W$.
The proof is finished.
\eo

The lemma below guarantees the convergence of solution $\Phi_{\epsilon}(t,\w,u_{\epsilon0})\ra\Phi_0(t,u_0)$ provided that $u_{\epsilon0}\ra u_0$ in $V$ as $\epsilon\ra0^+$.
This greatly helps prove the upper semi-continuity of random attractors.

\bl\label{le3.15} Let $u_\epsilon$ and $u$ be the solutions to the problem \eqref{1.1} -- \eqref{1.3} and \eqref{3.47} with initial conditions $u_{\epsilon0}$ and $u_0$, respectively. Then for $u_{\epsilon0}\ra u_0$ (as $\epsilon\ra0^+$) in $V$ and $t\geqslant0$, we have
\be\label{3.52}\lim_{\epsilon\ra0^+}u_\epsilon(t;\w,u_{\epsilon0})=u(t;u_0)\hs\mb{in }V,\hs\mb{for }\bP\mb{-a.s. }\w\in\W.\ee
\el
\bo
Let
\be\label{3.53}v_\epsilon(t;\w,v_{\epsilon0})=\me^{-\epsilon z(\theta_t\w)}u_{\epsilon}(t;\w,u_{\epsilon0})\,\mb{ and }\, \cY=v_\epsilon-u.\ee
Then by \eqref{2.7} and \eqref{3.47}, $\cY$ satisfies the following equation,
\be\label{3.54}
\frac{\pa \cY}{\pa t}+\De^2 \cY+2\De \cY+a\cY+\me^{2\epsilon z(\theta_t\w)}v_\epsilon G*v_{\epsilon}^2-uG*u^2+(\me^{-\epsilon z(\theta_t\w)}-1)h=\epsilon v_\epsilon z(\theta_t\w),
\ee
for $x\in U$ and $t>0$, with $\cY(0,\w)=v_{\epsilon0}(\w)-u_0=\me^{-\epsilon z(\w)}u_{\epsilon0}-u_0$.

Multiplying \eqref{3.54} by $\De^2\cY$ and integrating over $U$, we have
\be\label{3.55}\ba{c}\disp\frac12\frac{\di}{\di t}\|\De \cY\|^2+\|\De^2\cY\|^2+a\|\De \cY\|^2=-2(\De \cY,\De^2\cY)-(\me^{-\epsilon z(\theta_t\w)}-1)(h,\De^2 \cY)\\[1ex]
-(\me^{2\epsilon z(\theta_t\w)}v_\epsilon G*v_\epsilon^2-uG*u^2,\De^2\cY)+\epsilon z(\theta_t\w)(v_\epsilon,\De^2\cY).\ea
\ee
It is easy to see that
\be|2(\De \cY,\De^2 \cY)|\leqslant2\|\De \cY\|^2+\frac12\|\De^2 \cY\|^2,\ee
\be|\me^{-\epsilon z(\theta_t\w)}-1||(h,\De^2 \cY)|\leqslant\frac12(\me^{-\epsilon z(\theta_t\w)}-1)^2\|h\|^2+\frac12\|\De^2 \cY\|^2,\ee
\be|(v_\epsilon,\De^2\cY)|=|(\cY+u,\De^2 \cY)|
\leqslant\|\De \cY\|^2+|(\De u,\De \cY)|\leqslant\frac32\|\De \cY\|^2+\frac12\|\De u\|^2\hs\mb{and}
\ee
\be\ba{rl}&-(\me^{2\epsilon z(\theta_t\w)}v_\epsilon G*v_\epsilon^2-uG*u^2,\De^2\cY)\\
=&-\me^{2\epsilon z(\theta_t\w)}(\cY G*v_\epsilon^2,\De^2\cY)-\me^{2\epsilon z(\theta_t\w)}(uG*(v_\epsilon^2-u^2),\De^2\cY)\\
&-(\me^{2\epsilon z(\theta_t\w)}-1)(uG*u^2,\De^2\cY)\\
:=&I_1+I_2+I_3.\ea\ee
Now we estimate $I_1$, $I_2$ and $I_3$ respectively as follows.
By the assumption on $G$, Cauchy-Schwarz inequality, and the embedding inequality of Sobolev spaces, we have
\begin{align}
|I_1|=&\me^{2\epsilon z(\theta_t\w)}\left|(\De\cY G*v_\epsilon^2+2\na \cY\na G*v_\epsilon^2+\cY\De G*v_\epsilon^2,\De \cY)\right|\notag\\[1ex]
\leqslant&\b\me^{2\epsilon z(\theta_t\w)}\|v_\epsilon\|^2\int_U(|\De \cY|+2|\na \cY|+|\cY|)|\De \cY|\di x\notag\\
\lesssim&\me^{2\epsilon z(\theta_t\w)}\|v_\epsilon\|^2\|\De \cY\|^2.
\end{align}
Similarly, we also obtain
\begin{align}|I_2|&\leqslant\me^{2\epsilon z(\theta_t\w)}\(\De uG*(v_\epsilon^2-u^2)+2\na u\na G*(v_\epsilon^2-u^2)+u\De G*(v_\epsilon^2-u^2),\De \cY\)\notag\\
&\leqslant\b\me^{2\epsilon z(\theta_t\w)}\|v_\epsilon^2-u^2\|_{L^1(U)}\int_U(|\De u|+2|\na u|+|u|)|\De \cY|\di x\notag\\
&\lesssim\me^{2\epsilon z(\theta_t\w)}\|\cY\|\(\|v_\epsilon\|+\|u\|\)\|\De u\|\|\De \cY\|\notag\\
&\lesssim\me^{2\epsilon z(\theta_t\w)}\(\|v_\epsilon\|^2+\|\De u\|^2\)\|\De \cY\|^2
\end{align}
and
\begin{align}|I_3|
=&\left|1-\me^{2\epsilon z(\theta_t\w)}\right||(\De uG*u^2+2\na u\na G*u^2+u\De G*u^2,\De \cY)|\notag\\
\leqslant&\b\left|1-\me^{2\epsilon z(\theta_t\w)}\right|\|u\|^2\int_U(|\De u|+2|\na u|+|u|)|\De \cY|\di x\notag\\
\lesssim&\left|1-\me^{2\epsilon z(\theta_t\w)}\right|\|u\|^2\|\De u\|\|\De \cY\|\notag\\
\lesssim&\left|1-\me^{2\epsilon z(\theta_t\w)}\right|\|\De u\|^2\|u\|^2+\(1+\me^{2\epsilon z(\theta_t\w)}\)\|\De u\|^2\|\De \cY\|^2.\label{3.62}
\end{align}

Combining \eqref{3.55} -- \eqref{3.62}, we have the following inequality,
\begin{align}\frac{\di}{\di t}\|\De \cY\|^2
\leqslant&c\left[1+\epsilon|z(\theta_t\w)|+\me^{2\epsilon z(\theta_t\w)}\|v_\epsilon\|^2+\(1+\me^{2\epsilon z(\theta_t\w)}\)\|\De u\|^2\right]\|\De \cY\|^2\notag\\
&+c\left[(\me^{-\epsilon z(\theta_t\w)}-1)^2\|h\|^2+\left|\me^{2\epsilon z(\theta_t\w)}-1\right|\|\De u\|^2\|u\|^2+\epsilon|z(\theta_t\w)|\|\De u\|^2\right]\notag\\
:=&L_1(v_\epsilon,u,\epsilon z(\theta_t\w))\|\De \cY\|^2+L_2(u,\epsilon z(\theta_t\w)).\label{3.63}
\end{align}
Using Gronwall's lemma to \eqref{3.63}, we have for each fixed $t>0$,
\be\label{3.64}\ba{rl}
\|\De \cY(t)\|^2\leqslant&\me^{\int_0^tL_1(v_\epsilon(\vsig),u(\vsig),\epsilon z(\theta_\vsig\w))\di\vsig}\|\De \cY(0)\|^2\\
&\disp+\int_0^t\me^{\int_s^tL_1(v_\epsilon(\vsig),u(\vsig),\epsilon z(\theta_\vsig\w))\di\vsig}L_2(u(s),\epsilon z(\theta_s\w))\di s.\ea
\ee
Since $u_{\epsilon0}\ra u_0$ in $V$ as $\epsilon\ra0^+$, we have for $\bP$-a.s. $\w\in\W$, as $\epsilon\ra0$,
\begin{align}\|\De \cY(0)\|^2=&\|\De v_{\epsilon0}(\w)-\De u_0\|^2\notag\\
=&\|\me^{-\epsilon z(\w)}(\De u_{\epsilon0}-\De u_0)+(\me^{-\epsilon z(\w)}-1)\De u_0\|^2\notag\\
\leqslant&2\(\me^{-2\epsilon z(\w)}\|\De u_{\epsilon0}-\De u_0\|^2+|\me^{\epsilon z(\w)}-1|^2\|\De u_0\|^2\)\ra0,\label{3.65}\end{align}
as $\epsilon\ra0^+$.
Following the discussions in Subsection \ref{ss3.2}, we have the following estimate by \eqref{3.11} for the deterministic case,
\be\label{3.66}\|u(t)\|^2+\me^{-2t}\int_0^t\me^{2s}\|\De u(t)\|^2\di s\leqslant\me^{-2t}\|u_0\|^2+c(1+\|h\|^2).\ee
By \eqref{3.11} and assuming $0<\epsilon\leqslant1$, we have
\begin{align}&\int_0^t\me^{2\epsilon z(\theta_s\w)}\|v_\epsilon(s)\|^2+\(1+\me^{2\epsilon z(\theta_s\w)}\)\|\De u(s)\|^2\di s\notag\\
\leqslant&\me^{-2\epsilon z(\w)}\|u_{\epsilon0}\|^2\int_0^t\me^{-2s+2\epsilon\int_0^s z(\theta_{\vsig}\w)\di\vsig+2\epsilon z(\theta_s\w)}\di s\notag\\
&+c\(1+\|h\|^2\)\int_0^t\me^{-2s+2\epsilon z(\theta_s\w)}\int_0^s\me^{2r+2\epsilon\int_r^s z(\theta_{\vsig}\w)\di\vsig-2\epsilon z(\theta_{r}\w)}\di r\di s\notag\\
&+\(1+\me^{2\max_{0\leqslant s\leqslant t}|z(\theta_s\w)|}\)\int_0^t\me^{2s}\|\De u\|^2\di s.\label{3.67}
\end{align}
According to \eqref{3.66}, \eqref{3.67} and \eqref{2.5}, we can deduce that
\be\me^{\int_0^tL_1(v_\epsilon,u,\epsilon z(\theta_\vsig\w))\di\vsig}\leqslant R_0(t,\w)<\8\label{3.68}\ee
uniformly in $\epsilon$ for some positive $R_0(t,\w)$, since $\|u_{\epsilon0}\|$ is uniformly bounded in $\epsilon$ by the fact $u_{\epsilon0}\ra u_0$ in $V$.
Moreover, by \eqref{3.68}, we have
\begin{align*}&\int_0^t\me^{\int_s^tL_1(v_\epsilon,u,\epsilon z(\theta_\vsig\w))\di\vsig}L_2(u,\epsilon z(\theta_s\w))\di s\\
\leqslant&R_0(t,\w)\left[\max_{0\leqslant s\leqslant t}|\me^{-\epsilon z(\theta_s\w)}-1|^2t\|h\|^2+\epsilon\max_{0\leqslant s\leqslant t}|z(\theta_s\w)|\int_0^t\me^{2s}\|\De u(s)\|^2\di s\right]\\
&+R_0(t,\w)\max_{0\leqslant s\leqslant t}|\me^{-\epsilon z(\theta_s\w)}-1|\int_0^t\|\De u\|^2\|u\|^2\di s.
\end{align*}
By \eqref{3.66}, we have
$$\int_0^t\|\De u\|^2\|u\|^2\di s\leqslant\left[\|u_0\|^2+c\(1+\|h\|^2\)\right]\left[\|u_0\|^2+c\me^{2t}\(1+\|h\|^2\)\right]$$
and hence
\be\label{3.69}\int_0^t\me^{\int_s^tL_1(v_\epsilon,u,\epsilon z(\theta_\vsig\w))\di\vsig}L_2(u,\epsilon z(\theta_s\w))\di s\ra0\hs\mb{as }\epsilon\ra0^+.\ee

Combining \eqref{3.64}, \eqref{3.65}, \eqref{3.68} and \eqref{3.69} together, we obtain that
\be\label{3.70}\|\De v_\epsilon(t)-\De u(t)\|^2=\|\De \cY(t)\|^2\ra0,\hs\mb{as }\,\epsilon\ra0^+.\ee
Eventually, by \eqref{3.53} and \eqref{3.70}, we have
$$\ba{rl}\|\De u_\epsilon-\De u\|&=\|\me^{\epsilon z(\theta_t\w)}\De v_\epsilon-\De u\|=\|\me^{\epsilon z(\theta_t\w)}(\De v_\epsilon-\De u)+(\me^{\epsilon z(\theta_t\w)}-1)\De u\|\\
&\leqslant\me^{\epsilon z(\theta_t\w)}\|\De v_\epsilon-\De u\|+|\me^{\epsilon z(\theta_t\w)}-1|\|\De u\|\ra0,\ea
$$
as $\epsilon\ra0^+$, which immediately implies \eqref{3.52}.
The proof is complete.
\eo

By Proposition \ref{p3.6}, Lemmas \ref{le3.13}, \ref{le3.14} and \ref{le3.15}, we obtain the final conclusion.
\bt Assume that $h\in H$ and $G$ satisfies ({\bf Gp}).
Then for $\bP$-a.s. $\w\in\W$, we have
$$\lim_{\epsilon\ra0^+}{\rm dist}_{V}(\sA_\epsilon(\w),\sA_0)=0.$$
\et
\subsubsection{The case of non-negative kernel}

Now we consider the nonlocal stochastic SHE with the special non-negative kernel $G=J_{\varrho_0}$, which implies \eqref{2.12}.
Observe the similarity of the estimates in Lemmas \ref{le3.7}, \ref{le3.8} and \ref{le3.9}.
The results in Lemmas \ref{le3.13} and \ref{le3.14} trivially hold true for the problem \eqref{1.1} -- \eqref{1.3} when \eqref{2.12} holds for $G$ and $h\in H$.

In order to obtain the similar result in Lemma \ref{le3.15}, we follow the proof of Lemma \ref{le3.15}.
What are different lie in the estimations in \eqref{3.66} and \eqref{3.67}.
As a substitute of \eqref{3.66} for $u$, by \eqref{3.46}, we obtain
\begin{align}\|u(t)\|^2+\me^{-2t}\int_0^t\me^{2s}\|\De u(t)\|^2\di s\leqslant&\me^{-2t}\|u_0\|^2+c(1+\|h\|^2)+c\notag\\
\leqslant&\me^{-2t}\|u_0\|^2+c(1+\|h\|^2)\label{3.71}\end{align}
By \eqref{3.46} and $0<\epsilon\leqslant1$, we have a substitute of \eqref{3.67} as follows,
\begin{align}&
  \int_0^t\me^{2\epsilon z(\theta_s\w)}\|v_\epsilon(s)\|^2+\(1+\me^{2\epsilon z(\theta_s\w)}\)\|\De u(s)\|^2\di s\notag\\
\leqslant&
  \me^{-2\epsilon z(\w)}\|u_{\epsilon0}(\w)\|^2\int_0^t\me^{-2s+2\epsilon\int_0^s z(\theta_{\vsig}\w)\di\vsig+2\epsilon z(\theta_s\w)}\di s\notag\\
&
  +c\(1+\|h\|^2\)\int_0^t\me^{-2s+2\epsilon z(\theta_s\w)}\int_0^s\me^{2r+2\epsilon\int_r^s z(\theta_{\vsig}\w)\di\vsig-2\epsilon z(\theta_{r}\w)}\di r\di s\notag\\
&
  +c\int_0^t\me^{-2s+2\epsilon z(\theta_s\w)}\int_0^s\me^{2r+2\epsilon\int_r^s z(\theta_{\vsig}\w)\di\vsig+2\epsilon z(\theta_{r}\w)}\di r\di s\notag\\
&
  +\(1+\me^{2\max_{0\leqslant s\leqslant t}|z(\theta_s\w)|}\)\int_0^t\me^{2s}\|\De u\|^2\di s.\label{3.72}
\end{align}
The estimates \eqref{3.71}, \eqref{3.72} and \eqref{2.5} imply that \eqref{3.68} also holds for for this case when $G=J_{\varrho_0}$, with a different $R'_0(t,\w)$ instead of $R_0(t,\w)$.
Then using the same argument after \eqref{3.68} in the proof of Lemma \ref{le3.15}, one can smoothly prove for this case that, for $u_{\epsilon0}\ra u_0$ in $V$, as $\epsilon\ra0^+$ and for $t\geqslant0$, we have
$$\lim_{\epsilon\ra0^+}u_\epsilon(t;\w,u_{\epsilon0})=u(t;u_0)\hs\mb{in }V,\hs\mb{for }\bP\mb{-a.s. }\w\in\W.$$

The discussion above assures us the following upper semi-continuity of random attractors for the process generated by the problem \eqref{1.1} -- \eqref{1.3} with the special non-negative kernel $G=J_{\varrho_0}$.
\bt Assume that $h\in H$ and $G$ satisfies ({\bf Gn}).
Then for $\bP$-a.s. $\w\in\W$, we have
$$\lim_{\epsilon\ra0^+}{\rm dist}_{V}(\sA_\epsilon(\w),\sA_0)=0.$$
\et

\section{Existence of Ergodic Invariant Measures}
In this section, we consider the existence of ergodic invariant measures for the nonlocal stochastic SHE \eqref{1.1} in $H=L^2(U)$.
First we introduce some necessary notations and results.
\subsection{Preliminaries for ergodic invariant measures}
Let $X$ be a Banach space with its Borel $\sig$-algebra $\cB(X)$.
We use $\cB_{\rm b}(X)$ to denote the space of bounded Borel measurable functions on $X$.
In the following we rewrite the $X$-valued stochastic process $\cX(t,\w,x)$ as $\cX(t,x)(\w)$ for notational convenience with initial datum $x\in X$ and $\w\in\W$.
For a set $\Gam\in\cB(H)$, we define the {\em transition functions}
$$\cP_t(x,\Gam)=\bP(\cX(t,x)\in\Gam)\hs\mb{for all }t\geqslant0.$$
The {\em Markovian transition semigroup} $\cP_t$ on $\cB_{\rm b}(X)$ is defined as
$$\cP_t\vp(x)=\bE\vp(\cX(t,x))=\int_{X}\vp(y)\cP_t(x,\di y),\hs t\geqslant0,\;\vp\in\cB_{\rm b}(X),\;x\in X.$$
Note also that $\cP_t(x,\Gam)=\cP_t\1_{\Gam}(x)$.
The dual semigroup $\cP^*_t$ of $\cP_t$ is defined on and into the set of Borel probability measures $\nu$ on $X$ by
$$\cP^*_t\nu(\Gam)=\int_{X}\cP_t(x,\Gam)\nu(\di x),\hs\mb{for each }\Gam\in\cB(X).$$

The Markovian transition semigroup $\cP_t$ ($t\geqslant0$) is said to be {\em Feller} if for arbitrary $\vp\in \cC_{\rm b}(X)$
and $t\geqslant0$, the mapping $x\mapsto\cP_t\vp(x)$ is continuous.
An {\em invariant measure} for the stochastic process $\cX(t,x)$ is a probability measure $\nu$ on $X$, which is a fixed point for $\cP^*_t$, that is to say,
$$\int_X\cP_t(x,\Gam)\nu(\di x)=\nu(\Gam),\hs\mb{for each }\Gam\in\cB(X)\mb{ and }t\in\R^+.$$
Let $\nu$ be an invariant measure for $\cP_t$.
We say that $\nu$ is {\em ergodic} if
$$\lim_{T\ra\8}\frac1T\int_0^T\cP_t\vp\di t=\int_H\vp(x)\nu(\di x)\hs\mb{for all }\vp\in L^2(X,\mu).$$

For the discussion next, we need the following ``stochastic Gronwall lemma" (see \cite[Lemma 5.3]{GHZ09}).
\bp\label{p4.1} Fix $T>0$ and assume that $\cX$, $\cY$, $\cZ$, $\cR:[0,T)\X\W\ra\R$ are real valued, non-negative stochastic processes.
Let $\tau<T$ be a stopping time so that
$$\bE\int_0^\tau(\cR\cX+\cZ)\di s<\8.$$
Assume moreover that for some fixed deterministic constant $\kappa$,
$$\int_0^\tau\cR\di s<\kappa,\hs\mb{a.s.}$$
Suppose that
\be\label{4.1}\ba{l}\mb{ for all stopping times }0\leqslant\tau'<\tau''\leqslant\tau,\\
\disp\bE\(\sup_{s\in[\tau',\tau'']}\cX+\int_{\tau'}^{\tau''}\cY\di s\)\leqslant C_0\bE\(\cX(\tau')+\int_{\tau'}^{\tau''}(\cR\cX+\cZ)\di s\),\ea\ee
where $C_0$ is a positive constant independent of the choice of $\tau'$, $\tau''$.
Then
$$\bE\(\sup_{s\in[0,\tau]}\cX+\int_0^\tau\cY\di s\)\leqslant C\bE\(\cX(0)+\int_0^\tau\cZ\di s\),$$
where $C=C(C_0,T,\kappa)$ is a positive constant.
\ep

For the practical applications in the following, we sometimes need to relax the conditions given in Proposition \ref{p4.1}.
Actually, the relaxed conditions are equivalent to those in Proposition \ref{p4.1}, the details of which is stated in the following lemma.
\bl\label{le4.2} Under the conditions of Proposition \ref{p4.1}, the condition \eqref{4.1} is equivalent to the following condition:
\be\label{4.2}\ba{l}\mb{there is }\ve>0\mb{ such that for all stopping times }0\leqslant\tau'<\tau''\leqslant\tau\mb{ with }\tau''-\tau'<\ve,\\
\disp\bE\(\sup_{s\in[\tau',\tau'']}\cX+\int_{\tau'}^{\tau''}\cY\di s\)\leqslant C_\ve\bE\(\cX(\tau')+\int_{\tau'}^{\tau''}(\cR\cX+\cZ)\di s\),\ea\ee
where $C_\ve$ depends on $\ve$, but does not depend on the specific distributions of $\tau'$ and $\tau''$.\el

\bo Obviously, the condition \eqref{4.1} implies \eqref{4.2} by arbitrary choice of $\ve>0$.
Now we show the converse.

Let \eqref{4.2} hold true.
Select a finite sequence of stopping times
$$0=\tilde\tau_0<\tilde\tau_1<\cdots<\tilde\tau_N<\tilde\tau_{N+1}=T$$
such that $\tilde\tau_{k}-\tilde\tau_{k-1}<\ve$.
Pick arbitrarily $\tau'$ and $\tau''$ with $0\leqslant\tau'<\tau''\leqslant\tau$.
Surely, we have $0\leqslant k_1\leqslant k_2\leqslant N$ such that
$$\tau'\in[\tilde\tau_{k_1},\tilde\tau_{k_1+1})\cap[0,\tau]
\hs
\mb{and}
\hs
\tau''\in(\tilde\tau_{k_2},\tilde\tau_{k_2+1}]\cap[0,\tau].$$
We claim that
\be\bE\(\sup_{s\in[\tau',\tau'']}\cX+\int_{\tau'}^{\tau''}\cY\di s\)\leqslant \(1+C_\ve\)^{k_2-k_1+1}\bE\(\cX(\tau')+\int_{\tau'}^{\tau''}(\cR\cX+\cZ)\di s\),\label{4.3}\ee
which implies \eqref{4.1} immediately.

We show the claim by induction.
Firstly, by \eqref{4.2}, since $\tilde\tau_{k_1+1}-\tau'<\ve$, we have
$$\bE\(\sup_{s\in[\tau',\tilde\tau_{k_1+1}]}\cX+\int_{\tau'}^{\tilde\tau_{k_1+1}}\cY\di s\)\leqslant \(1+C_\ve\)\bE\(\cX(\tau')+\int_{\tau'}^{\tilde\tau_{k_1+1}}(\cR\cX+\cZ)\di s\).$$
Assume that for $i\in[1,k_2-k_1]$,
$$\bE\(\sup_{s\in[\tau',\tilde\tau_{k_1+i}]}\cX+\int_{\tau'}^{\tilde\tau_{k_1+i}}\cY\di s\)\leqslant \(1+C_\ve\)^i\bE\(\cX(\tau')+\int_{\tau'}^{\tilde\tau_{k_1+i}}(\cR\cX+\cZ)\di s\).$$
Then by \eqref{4.2},
\begin{align*}&\bE\(\sup_{s\in[\tau',\tilde\tau_{k_1+i+1}\wedge\tau'']}\cX+\int_{\tau'}^{\tilde\tau_{k_1+i+1}\wedge\tau''}\cY\di s\)\\
\leqslant&
  \(1+C_\ve\)^i\bE\(\cX(\tau')+\int_{\tau'}^{\tilde\tau_{k_1+i}}(\cR\cX+\cZ)\di s\)+\bE\(\sup_{s\in[\tilde\tau_{k_1+i},\tilde\tau_{k_1+i+1}\wedge\tau'']}\cX+\int_{\tilde\tau_{k_1+i}}^{\tilde\tau_{k_1+i+1}\wedge\tau''}\cY\di s\)\\
\leqslant&
  \(1+C_\ve\)^i\bE\(\cX(\tau')+\int_{\tau'}^{\tilde\tau_{k_1+i}}(\cR\cX+\cZ)\di s\)\!+\!C_\ve\bE\(\cX(\tilde\tau_{k_1+i})+\int_{\tilde\tau_{k_1+i}}^{\tilde\tau_{k_1+i+1}\wedge\tau''}(\cR\cX+\cZ)\di s\)\\
\leqslant&
  \(1+C_\ve\)^{i+1}\bE\(\cX(\tau')+\int_{\tau'}^{\tilde\tau_{k_1+i}}(\cR\cX+\cZ)\di s\)+C_\ve\bE\int_{\tilde\tau_{k_1+i}}^{\tilde\tau_{k_1+i+1}\wedge\tau''}(\cR\cX+\cZ)\di s\\
\leqslant&
  \(1+C_\ve\)^{i+1}\bE\(\cX(\tau')+\int_{\tau'}^{\tilde\tau_{k_1+i+1}\wedge\tau''}(\cR\cX+\cZ)\di s\),
\end{align*}
which completes the proof of the claim.

Note that $k_2-k_1\leqslant N$.
Hence \eqref{4.1} follows from \eqref{4.3} at once.
\eo
\Vs

For the following arrangement, we first consider the nonlocal stochastic Swift-Hohenberg model \eqref{1.1} -- \eqref{1.3} with a positive kernel, saying, under the assumption \eqref{2.10}.
For the case when the kernel $G$ is non-negative, we again only consider the special non-negative kernel $G=J_{\varrho_0}$, which hence satisfies \eqref{2.12}.
For this case, we will give a succinct discussion afterwards.

\subsection{The case of positive kernel}

Now we get back to the original problem and consider the stochastic process on $H$ given by \eqref{3.1} generated by the problem \eqref{1.1} -- \eqref{1.3} under the assumption ({\bf Gp}).
Here we always assume that $\epsilon\in(0,1]$ and $h\in H$.
We can infer from \eqref{1.1} that the solution $u(t,u_0)(\w)=u(t;\w,u_0)$ satisfies the following integral equation, for all $t\geqslant s\geqslant0$,
\be\label{4.A}u(t)=\me^{-A(t-s)}u(s)-\int_s^t\me^{-A(t-\vsig)}[2A^{\frac12}u+au+uG*u^2+h]\di\vsig+\epsilon\int_s^t\me^{-A(t-\vsig)}u\di W_\vsig.\ee
Moreover, note that the Wiener process $W_t$ has a natural filtration
$$\cF_t:=\vsig(\w(s):0\leqslant s\leqslant t),\hs\mb{for all $t\geqslant 0$}.$$
By the representation \eqref{2.13} of the solution, we know that for each deterministic initial datum $u_0\in H$,
the solution $u(t,u_0)$ is an $H$-valued predictable process with respect to the filtration $\~\cF_t$ with $\~\cF_0=\{\es,\W\}$ and $\~\cF_t=\cF_t$ for $t>0$.

In order to prove the existence of invariant measures, we need some auxiliary consequences of the estimates of moment bounds and some related probabilities for the solution $u(t,u_0)$.

\subsubsection{Estimates on moment bounds of solutions}

In this part, we allow the initial data $u_0$, $u_{10}$ and $u_{20}$ to be $H$-valued random variables.

\bl\label{le4.3} For each solution $u(t,u_0)$ of \eqref{1.1} -- \eqref{1.3}, $p\geqslant1$ and $t>0$ with $u_0\in L^{2p}(\W,H)$, we have
\be\label{4.4}
\bE\(\sup_{s\in[0,t]}\|u(s)\|^{2p}+\int_0^t\|u\|^{2p-2}\|\De u\|^2\di s\)\lesssim\bE\|u_0\|^{2p}+(\|h\|^{2p}+1)t.\ee
\el
\bo Let $g(u):=\|u\|^{2p}$ ($p\geqslant1$) for each solution $u(t,u_0)$.
Note that the derivatives of $g(u)$ satisfies
\[
g_u=2p\|u\|^{2p-2}(u,\cdot)\hs\mb{and}\hs g_{uu}=2p(2p-2)\|u\|^{2p-4}(u,\cdot)^2+2p\|u\|^{2p-2}(\cdot,\cdot).
\]
Applying the It\^o's Formula to the original equation \eqref{1.1}, we have (replacing the time $t$ by $\vsig$)
\begin{align}
&\di\|u\|^{2p}+\|u\|^{2p-2}\{2p[\|\De u\|^2+2(\De u,u)+a\|u\|^2+(u^2,G*u^2)+(h,u)]-p(2p-1)\epsilon^2\|u\|^2\}\di\vsig\notag\\
=&2p\epsilon\|u\|^{2p}\di W_\vsig,\label{4.5}
\end{align}
which implies, by integration over the time interval $[0,s]$ and inequalities similar to \eqref{3.9} and \eqref{3.10}, that
\begin{align}&\|u(s)\|^{2p}+\int_0^s\|u\|^{2p-2}\|\De u\|^2\di\vsig\notag\\
\leqslant
&\|u_0\|^{2p}+\int_0^s[c(\|u\|^{2p}+\|h\|\|u\|^{2p-1})-\al\|u\|^{2p+2}]\di\vsig+2p\epsilon\int_0^s\|u\|^{2p}\di W_\vsig\notag\\
\leqslant&\|u_0\|^{2p}+\int_0^s[\|h\|^{2p}+c\|u\|^{2p}-\al\|u\|^{2p+2}]\di\vsig+2p\epsilon\int_0^s\|u\|^{2p}\di W_\vsig.\label{4.6}
\end{align}
By considering the supremum over $[0,t]$, we obtain
\begin{align*}\sup_{s\in[0,t]}\|u(s)\|^{2p}+\int_0^t\|u\|^{2p-2}\|\De u\|^2\di s\leqslant&\|u_0\|^{2p}+\int_0^t[\|h\|^{2p}+c\|u\|^{2p}-\al\|u\|^{2p+2}]\di s\\
&+2p\epsilon\sup_{s\in[0,t]}\left|\int_0^s\|u\|^{2p}\di W_s\right|.\end{align*}
Taking the expected value and using the Burkholder-Davis-Gundy inequality, we have
\begin{align*}&\bE\(\sup_{s\in[0,t]}\|u(s)\|^{2p}+\int_0^t\|u\|^{2p-2}\|\De u\|^2\di s\)\\
\leqslant&
  \bE\|u_0\|^{2p}+\int_0^t[\|h\|^{2p}+c\bE\|u\|^{2p}-\al\bE\|u\|^{2p+2}]\di s+2p\epsilon\bE\(\int_0^t\|u\|^{4p}\di s\)^{\frac12}\\
\leqslant&
  \bE\|u_0\|^{2p}+\int_0^t[\|h\|^{2p}+c\bE\|u\|^{2p}-\al\bE\|u\|^{2p+2}]\di s+\frac12\bE\sup_{s\in[0,t]}\|u(s)\|^{2p}+c\int_0^t\bE\|u(s)\|^{2p}\di s,
\end{align*}
which means
\begin{align}&\bE\(\sup_{s\in[0,t]}\|u(s)\|^{2p}+\int_0^t\|u\|^{2p-2}\|\De u\|^2\di s\)\notag\\
\lesssim&\bE\|u_0\|^{2p}+\int_0^t[\|h\|^{2p}+\bE\(c\|u\|^{2p}-\al\|u\|^{2p+2}\)]\di s.\label{4.7}\end{align}
Since $c\|u\|^{2p}-\al\|u\|^{2p+2}$ has an upper bound that is independent of $\|u\|$, the conclusion \eqref{4.4} follows immediately from \eqref{4.7}.
The proof is finished.
\eo

Now we define a stopping time $\tau^p_\kappa$ for $p\geqslant1$, $\kappa>0$ and $u_{0}\in L^{2p}(\W,H)$,
$$\tau^p_\kappa(u_{0}):=\inf\left\{t\geqslant0:\int_0^t\|u(s,u_{0})\|^{2p}\di s>\kappa\right\},$$
where $u(t,u_{0})$ denotes the strong solution to \eqref{1.1} -- \eqref{1.3} with the initial datum $u_{0}$.
And for arbitrary $u_{10}$, $u_{20}\in H$, we denote
\be\label{4.8}\tau^p_\kappa(u_{10},u_{20}):=\tau^p_\kappa(u_{10})\wedge\tau^p_\kappa(u_{20}),\ee
where $a\wedge b$ denote the smaller one of $a$ and $b$ for all $a,\,b\in\R$.

\bl\label{le4.4} For all $u_0\in L^{2p}(\W,H)$, $\kappa>0$, $p\geqslant1$ and $t>0$, we have
\be\label{4.9}\bP(\tau_{\kappa}^p(u_0)<t)\lesssim\kappa^{-1}t\left[\bE\|u_0\|^{2p}+(\|h\|^{2p}+1)t\right].\ee
\el

\bo Given each $u_0\in L^{2p}(\W,H)$, $\kappa>0$, $p\geqslant1$ and $t>0$, by Markovian inequality, we see that
$$
\bP(\tau_{\kappa}^p(u_0)<t)
=\bP\(\int_0^t\|u(s)\|^{2p}\di s>\kappa\)
\leqslant\bP\(\sup_{s\in[0,t]}\|u(s)\|^{2p}>\frac{\kappa}{t}\)
\leqslant\kappa^{-1}t\bE\sup_{s\in[0,t]}\|u(s)\|^{2p}.
$$
Thus \eqref{4.9} is a trivial deduction of the above estimate and Lemma \ref{le4.3}.
The proof is finished.
\eo

In order to obtain an apparent description of the instant parabolic regularization for the solutions to \eqref{1.1} -- \eqref{1.3},
we define another stopping time $\hat\tau_{\hat\kappa}(u_0)$ for each $u_0\in L^6(\W,H)$ and $\hat\kappa>0$ such that
$$\hat\tau_{\hat\kappa}(u_0)=\inf\{t\geqslant 0:t\|\De u(t)\|^2>\hat\kappa\},$$
and
\be\label{4.10}\hat\tau_{\hat\kappa}(u_{10},u_{20})=\hat\tau_{\hat\kappa}(u_{10})\wedge\hat\tau_{\hat\kappa}(u_{20}),\hs\mb{for all }u_{10},\,u_{20}\in H.\ee

To estimate the probability involving $\hat\tau_{\hat\kappa}$, we first give a suitable local moment bound on the $V$-norm of the solutions.

\bl\label{le4.5}For each $T>0$, $\kappa>0$ and $u_0\in L^6(\W,H)$, it holds that
\be\label{4.11}\bE\sup_{s\in[0,t\wedge\tau_{\kappa}^2(u_0)]}\(s\|\De u(s)\|^2\)\leqslant C_1\left[\bE\|u_0\|^2+(\|h\|^2+1)t+\|h\|^2t^2\right],\hs\mb{for all }t\in[0,T],\ee
where $C_1=C_1(T,\kappa)$ is a positive constant.
\el

\bo From the original equation \eqref{1.1}, we easily have for $s>0$,
$$\di(s^{\frac12}u)=[-A(s^{\frac12}u)+2^{-1}s^{-\frac12}u-s^{\frac12}(2\De u+au+uG*u^2+h)]\di s+\epsilon s^{\frac12}u\di W_s,$$
which indicates that the term $s^{\frac12}u$ satisfies for all stopping times $0\leqslant\tau'< s\leqslant\tau''\leqslant t\wedge\tau_\kappa^2(u_0)$,
\begin{align} s^{\frac12}u(s)=&\me^{-A(s-\tau')}\tau'^{\frac12}u(\tau')+\int_{\tau'}^s\me^{-A(s-\vsig)}[2^{-1}\vsig^{-\frac12}u-\vsig^{\frac12}(2\De u+au+uG*u^2+h)]\di\vsig\notag\\
&+\epsilon\int_{\tau'}^s\me^{-A(s-\vsig)}\vsig^{\frac12}u(\vsig)\di W_\vsig.\label{4.12}
\end{align}
Applying $-\De$ to \eqref{4.12}, estimating the $V$-norm and by \eqref{2.1}, \eqref{2.2}, we have
\begin{align}s^{\frac12}\|\De u(s)\|\leqslant&
  \|A^{\frac12}\me^{-A(s-\tau')}\tau'^{\frac12}u(\tau')\|+\int_{\tau'}^{s}\left\|A^{\frac12}\me^{-A(s-\vsig)}\vsig^{\frac12}(A^{\frac12}u+au+uG*u^2+h)\right\|\di\vsig\notag\\
&
  +2^{-1}\int_{\tau'}^{t}\left\|A^{\frac12}\me^{-A(s-\vsig)}\vsig^{-\frac12}u\right\|\di\vsig+\epsilon\left\|\int_{\tau'}^{s}A^{\frac12}\me^{-A(s-\vsig)}\vsig^{\frac12}u(\vsig)\di W_\vsig\right\|\notag\\
\lesssim&
  \tau'^{\frac12}\|\De u(\tau')\|+\int_{\tau'}^{s}(s-\vsig)^{-\frac12}\vsig^{\frac12}(\|\De u\|+\|h\|)\di\vsig+\int_{\tau'}^{s}\vsig^{\frac12}\|\De u\|(1+\|u\|^2)\di\vsig\notag\\
&
+\int_{\tau'}^{s}(s-\vsig)^{-\frac12}\vsig^{-\frac12}\|u\|\di\vsig+\epsilon\left\|\int_{\tau'}^{s}\me^{-A(s-\vsig)}\vsig^{\frac12}\De u(\vsig)\di W_\vsig\right\|.\label{4.13}
\end{align}
Note that
$$\int_{\tau'}^{s}(s-\vsig)^{-\frac12}\vsig^{\frac12}\|\De u\|\di\vsig\leqslant2(s-\tau')^{\frac12}\(\sup_{\vsig\in[\tau',s]}\(\vsig\|\De u(\vsig)\|^2\)\)^{\frac12}\hs\mb{and}$$
\begin{align*}&\bE\sup_{s\in[\tau',\tau'']}\left\|\int_{\tau'}^{s}\me^{-A(s-\vsig)}\vsig^{\frac12}\De u(\vsig)\di W_\vsig\right\|^2
\lesssim\bE\int_{\tau'}^{\tau''}\left\|\me^{-A(s-\vsig)}\vsig^{\frac12}\De u(\vsig)\right\|^2\di\vsig\\
\lesssim&\bE\int_{\tau'}^{\tau''}\vsig\|\De u(\vsig)\|^2\di\vsig\lesssim(\tau''-\tau')\bE\sup_{\vsig\in[\tau',\tau'']}\(\vsig\|\De u(\vsig)\|^2\),\end{align*}
by Burkholder-Davis-Gundy inequality.
Thus, by considering the square of the two sides of \eqref{4.13} and the supremum over $t\in[\tau',\tau'']$, and then taking the expected value, we get
\begin{align}&\bE\sup_{s\in[\tau',\tau'']}\(s\|\De u(s)\|^2\)\notag\\
\lesssim&
  \bE\(\tau'\|\De u(\tau')\|^2\)+(\tau''-\tau')\bE\sup_{s\in[\tau',\tau'']}\(s\|\De u(s)\|^2\)+(\tau''-\tau')\bE\int_{\tau'}^{\tau''}\vsig\|\De u\|^2\(1+\|u\|^4\)\di\vsig\notag\\
&
  +\bE\sup_{s\in[\tau',\tau'']}\|u(s)\|^2\int_{\tau'}^{\tau''}(\tau''-\vsig)^{-\frac12}\vsig^{-\frac12}\di\vsig+\|h\|^2\int_{\tau'}^{\tau''}(\tau''-\vsig)^{-\frac12}\vsig^{\frac12}\di\vsig.\label{4.14}
\end{align}
By the definition of the denotation $\lesssim$, we can easily find a positive constant $\ve$ such that, whenever $\tau''-\tau'<\ve$, it follows from \eqref{4.14} that
\begin{align*}&\bE\sup_{t\in[\tau',\tau'']}\(t\|\De u(t)\|^2\)\\
\leqslant&\frac12\bE\sup_{t\in[\tau',\tau'']}\(t\|\De u(t)\|^2\)+
  c\bE\(\tau'\|\De u(\tau')\|^2\)+c\ve\bE\int_{\tau'}^{\tau''}\vsig\|\De u\|^2\(1+\|u\|^4\)\di\vsig\\
&
  +c\bE\int_{\tau'}^{\tau''}(\tau''-\vsig)^{-\frac12}\left[\vsig^{-\frac12}\sup_{s\in[\tau',\tau'']}\|u(s)\|^2+\vsig^{\frac12}\|h\|^2\right]\di\vsig,
\end{align*}
and hence
\begin{align*}\bE\sup_{t\in[\tau',\tau'']}\(t\|\De u(t)\|^2\)\lesssim&
  \bE\(\tau'\|\De u(\tau')\|^2\)+\ve\bE\int_{\tau'}^{\tau''}\vsig\|\De u\|^2\(1+\|u\|^4\)\di\vsig\\
&
  +\bE\int_{\tau'}^{\tau''}(\tau''-\vsig)^{-\frac12}\left[\vsig^{-\frac12}\sup_{s\in[\tau',\tau'']}\|u(s)\|^2+\vsig^{\frac12}\|h\|^2\right]\di\vsig.
\end{align*}
At last, by Proposition \ref{p4.1}, Lemmas \ref{le4.3} and \ref{le4.2}, we finally obtain \eqref{4.11},
where we have used the integrals
$$\int_0^t\(\frac{\vsig}{t-\vsig}\)^{\frac12}\di\vsig=\frac{\pi t}2 \hs\mb{and}\hs\int_0^t(t-\vsig)^{-\frac12}\vsig^{-\frac12}\di\vsig=\pi.$$
The proof is accomplished.\eo

\bl\label{le4.6} For each $u_0\in L^6(\W,H)$, $\hat\kappa,\kappa>0$ and $t>0$, we have
\be\label{4.15}\bP(\hat\tau_{\hat\kappa}(u_0)<t)\leqslant\hat\kappa^{-1}C_1(t,\kappa)[\bE\|u_0\|^2+(\|h\|^2+1)t+\|h\|^2t^2]+\kappa^{-1}ct[\bE\|u_0\|^4+(\|h\|^4+1)t],\ee
where $C_1$ is given in Lemma \ref{le4.5}.
\el

\bo By applying Markovian inequality, we have
\begin{align*}\bP(\hat\tau_{\hat\kappa}(u_0)<t)=&\bP\(\sup_{s\in[0,t]}(s\|\De u(s)\|^2)>\hat\kappa\)\\
\leqslant&\bP\(\sup_{s\in[0,t\wedge\tau_\kappa^2(u_0)]}(s\|\De u(s)\|^2)>\hat\kappa\)+\bP(\tau_\kappa^2(u_0)<t)\\
\leqslant&\hat\kappa^{-1}\bE\(\sup_{s\in[0,t\wedge\tau_\kappa^2(u_0)]}(s\|\De u(s)\|^2)\)+\bP(\tau_\kappa^2(u_0)<t),
\end{align*}
which infers \eqref{4.15} by Lemma \ref{le4.5} and Lemma \ref{le4.4} and ends the proof.
\eo
\subsubsection{Feller property}

\bl\label{le4.7} For each fixed deterministic $T>0$, $\kappa>0$ and $u_{10},\,u_{20}\in H$ (deterministic), it holds that
\be\label{4.16}\bE\sup_{s\in[0,t\wedge\tau_\kappa^1(u_{10},u_{20})]}\|u(s,u_{10})-u(s,u_{20})\|^2\leqslant C_2\|u_{10}-u_{20}\|^2,\hs\mb{for all }t\in[0,T],\ee
where $C_2=C_2(T,\kappa)$ is a positive constant.
\el
\bo
Let $t\in[0,T]$, $u_i(t):=u(t,u_{i0})$, $i=1,2$ and $\cZ(t)=u_1(t)-u_2(t)$, for $t\in[0,T]$.
Then by \eqref{1.1} we have
\[\di\cZ+[\De^2\cZ+2\De\cZ+a\cZ+\cZ G*u^2_2+u_1G*(u_2^2-u_1^2)]\di t=\epsilon \cZ\di W_t.\]
Applying It\^o's Formula to $\|\cZ(t)\|^2$ and similar to the last proof, we have for all stopping times $\tau'$ and $\tau''$ with $0\leqslant\tau'\leqslant s\leqslant\tau''\leqslant t\wedge\tau_\kappa^1(u_{10},u_{20})$,
\begin{align}
\|\cZ(s)\|^2=&
  \|\cZ(\tau')\|^2+\int_{\tau'}^{s}[\epsilon^2\|\cZ\|^2-2\|\De\cZ\|^2-4(\De\cZ,\cZ)-2a\|\cZ\|^2-2(\cZ^2,G*u_2^2)]\di s\notag\\
&
  -2\int_{\tau'}^{s}(\cZ u_1,G*(u_2^2-u_1^2))\di\vsig+2\epsilon\int_{\tau'}^{t}\|\cZ(\vsig)\|^2\di W_\vsig\notag\\
\leqslant&
  \|\cZ(\tau')\|^2\!\!+\suo\int_{\tau'}^{s}\![(\epsilon^2+2-2a)\|\cZ(\vsig)\|^2\!+\!2\beta\|\cZ(\vsig)\|\|u_1\|\|u_2^2-u_1^2\|]\di\vsig+2\epsilon\!\!\int_{\tau'}^{s}\|\cZ(\vsig)\|^2\di W_\vsig\notag\\
\leqslant&
  \|\cZ(\tau')\|^2+c\!\int_{\tau'}^{s}[\|\cZ\|^2+\|\cZ\|^2(\|u_1(\vsig)\|^2+\|u_2(\vsig)\|^2)]\di\vsig+2\epsilon\!\int_{\tau'}^{s}\|\cZ(\vsig)\|^2\di W_\vsig,\label{4.17}
\end{align}
where Young's inequality is used.
Taking a supremum over $t\in[\tau',\tau'']$ for \eqref{4.17} and then considering the expectation, we obtain
\begin{align}
\bE\sup_{s\in[\tau',\tau'']}\|\cZ(s)\|^2
\leqslant&
  \bE\|\cZ(\tau')\|^2+c\bE\int_{\tau'}^{\tau''}\|\cZ(\vsig)\|^2(\|u_1(\vsig)\|^2+\|u_2(\vsig)\|^2+1)\di\vsig\notag\\
&
  +2\bE\sup_{s\in[\tau',\tau'']}\left|\int_{\tau'}^s\|\cZ(\vsig)\|^2\di W_\vsig\right|\label{4.18}
\end{align}
Observe by the Burkholder-Davis-Gundy inequality that
\begin{align}
\bE\sup_{s\in[\tau',\tau'']}\left|\int_{\tau'}^s\|\cZ(\vsig)\|^2\di\vsig\right|\leqslant
&
  c\bE\(\int_{\tau'}^{\tau''}\|\cZ(\vsig)\|^4\di\vsig\)^{\frac12}\notag\\
\leqslant&
  c\bE\(\sup_{t\in[\tau',\tau'']}\|\cZ(t)\|^2\int_{\tau'}^{\tau''}\|\cZ(\vsig)\|^2\di\vsig\)^{\frac12}\notag\\
\leqslant&
  \frac14\bE\sup_{t\in[\tau',\tau'']}\|\cZ(t)\|^2+c\bE\int_{\tau'}^{\tau''}\|\cZ(\vsig)\|^2\di\vsig.\label{4.19}
\end{align}
The inequalities \eqref{4.18} and \eqref{4.19} indicate that
$$\bE\sup_{s\in[\tau',\tau'']}\|\cZ(s)\|^2
\lesssim\bE\|\cZ(\tau')\|^2+\bE\int_{\tau'}^{\tau''}\|\cZ(\vsig)\|^2(\|u_1(\vsig)\|^2+\|u_2(\vsig)\|^2+1)\di\vsig.$$
By \eqref{4.4} and the definition of $\tau_\kappa^1(u_{10},u_{20})$, we know
$$\bE\int_0^T\|\cZ(\vsig)\|^2(\|u_1(\vsig)\|^2+\|u_2(\vsig)\|^2+1)\di\vsig\lesssim [\|u_{10}\|^4+\|u_{20}\|^4+(\|h\|^4+1)T]T<\8,$$
$$\int_0^{t\wedge \tau_\kappa^1(u_{10},u_{20})}(\|u_1(\vsig)\|^2+\|u_2(\vsig)\|^2+1)\di\vsig\leqslant2\kappa+T.$$
Then by Proposition \ref{p4.1}, we obtain \eqref{4.16} and finish the proof.
\eo

Now we are well prepared to show the Feller property of the semigroup $\cP_t$, with the details stated in the next theorem.

\bt\label{th4.8} For each $\vp\in \cC_{\rm b}(H)$ and $t\geqslant0$, the mapping $u\mapsto\cP_t\vp(u)$ is continuous.
\et
\bo Given arbitrary $u_{10},\,u_{20}\in H$ and $\kappa_1,\hat\kappa>0$, we estimate $|\cP_t\vp(u_{10})-\cP_t\vp(u_{20})|$ as follows.
\begin{align}&|\cP_t\vp(u_{10})-\cP_t\vp(u_{20})|=|\bE(\vp(u(t,u_{10}))-\vp(u(t,u_{20})))|\notag\\
\leqslant&
  |\bE(\vp(u(t,u_{10}))-\vp(u(t,u_{20})))\1_{\{\tau_{\kappa_1}^1(u_{10},u_{20})<t\}}|\notag\\
&
  +|\bE(\vp(u(t,u_{10}))-\vp(u(t,u_{20})))\1_{\{\hat\tau_{\hat\kappa}(u_{10},u_{20})<t\}}|\notag\\
&
  +|\bE(\vp(u(t,u_{10}))-\vp(u(t,u_{20})))\1_{\{\tau_{\kappa_1}^1(u_{10},u_{20})\geqslant t\}}\1_{\{\hat\tau_{\hat\kappa}(u_{10},u_{20})\geqslant t\}}|\notag\\
:=&P_1+P_2+P_3,\label{4.20}
\end{align}
where $\tau_{\kappa_1}^1$ and $\hat\tau_{\hat\kappa}$ are defined as \eqref{4.8} and \eqref{4.10}, respectively.

Let $\|\vp\|_\8:=\sup_{u\in H}|\vp(u)|$ be the sup-norm of $\cC_{\rm b}(H)$ and the distance of the initial data satisfy $\|u_{10}-u_{20}\|<1$.
We first consider $P_1$.
By Lemma \ref{le4.4}, we have
\begin{align}P_1\leqslant&2\|\vp\|_{\8}\bP(\tau_{\kappa_1}^1(u_{10},u_{20})<t)\leqslant2\|\vp\|_{\8}\(\bP(\tau_{\kappa_1}^1(u_{10})<t)+\bP(\tau_{\kappa_1}^1(u_{20})<t)\)\notag\\
\lesssim&
  \kappa_1^{-1}\|\vp\|_\8 t[1+\|u_{10}\|^2+(\|h\|^2+1)t].\label{4.21}
\end{align}
For $P_2$, by Lemma \ref{le4.6}, we have for each $\kappa_2>0$,
\begin{align}P_2\leqslant&2\|\vp\|_{\8}\bP(\hat\tau_{\hat\kappa}(u_{10},u_{20})<t)\leqslant2\|\vp\|_{\8}\(\bP(\hat\tau_{\hat\kappa}(u_{10})<t)+\bP(\hat\tau_{\hat\kappa}(u_{20})<t)\)\notag\\
\lesssim&\hat\kappa^{-1}C_1(t,\kappa_2)\|\vp\|_\8[1+\|u_{10}\|^2+(\|h\|^2+1)t+\|h\|^2t^2]\notag\\
&
  +\kappa_2^{-1}\|\vp\|_\8 t[1+\|u_{10}\|^4+(\|h\|^4+1)t].\label{4.22}
\end{align}
For arbitrarily small $\ve>0$, by taking $\kappa_1=\kappa_2$ sufficiently large and fixed and then $\hat\kappa$ sufficiently large and fixed in \eqref{4.21} and \eqref{4.22}, we can ensure that
\be P_1+P_2<\frac\ve2.\ee

Now that we have fixed $\tau_{\kappa_1}^1$ and $\hat\kappa$, next we can delve into $P_3$.
For this, we approximate $\vp$ given above by a Lipschitz function $\tilde\vp$ (to be chosen below).
Observe that on the set $\{\hat\tau_{\hat\kappa}(u_{10},u_{20})\geqslant t\}$, we have $\|\De u(t,u_{i0})\|^2\leqslant\hat\kappa/t$.
As a result, $P_3$ can be bounded as follows,
\begin{align}P_3\leqslant&2\sup_{u\in\ol{\bf B}_V\((\hat\kappa/t)^{\frac12}\)}\|\vp(u)-\tilde\vp(u)\|+|\bE(\tilde\vp(u(t,u_{10}))-\tilde\vp(u(t,u_{20})))\1_{\{\tau_{\kappa_1}^1(u_{10},u_{20})\geqslant t\}}|\notag\\
\leqslant&
  2\suo\sup_{u\in\ol{\bf B}_V\((\hat\kappa/t)^{\frac12}\)}\suo\|\vp(u)-\tilde\vp(u)\|+L_{\tilde\vp}\bE\|\(u(t,u_{10})-u(t,u_{20})\)\1_{\{\tau_{\kappa_1}^1(u_{10},u_{20})\geqslant t\}}\|,
\end{align}
where $\ol{\bf B}_V(r)$ is the closed ball in $V$ centered at $0$ with radius $r$, and $L_{\tilde\vp}$ is the Lipschitz constant of $\tilde\vp$.

By the compact embedding of $V$ into $H$, we know the closed ball
$\ol{\bf B}_V\((\hat\kappa/t)^{\frac12}\)$ is compact in $H$.
By the density of Lipschitz functions in $\cC_{\rm b}\(\ol{\bf B}_V\((\hat\kappa/t)^{\frac12}\)\)$, we can find a Lipschitz function $\tilde\vp\in \cC_{\rm b}\(\ol{\bf B}_V\((\hat\kappa/t)^{\frac12}\)\)$ such that
\be\sup_{u\in\ol{\bf B}_V\((\hat\kappa/t)^{\frac12}\)}\|\vp(u)-\tilde\vp(u)\|<\frac\ve8.\ee
Now the choice of $\tilde\vp$ fixes $L_{\tilde\vp}$.
Then by Jensen's inequality and Lemma \ref{le4.7}, we have
\begin{align}&L_{\tilde\vp}\bE\|\(u(t,u_{10})-u(t,u_{20})\)\1_{\{\tau_{\kappa_1}^1(u_{10},u_{20})\geqslant t\}}\|\notag\\
\leqslant&
  L_{\tilde\vp}\(\bE\sup_{s\in[0,t\wedge\tau_{\kappa_1}^1(u_{10},u_{20})]}\|u(t,u_{10})-u(t,u_{20})\|^2\)^{\frac12}\notag\\
\leqslant&C_2(t,\kappa_1)L_{\tilde\vp}\|u_{10}-u_{20}\|.
\end{align}
Hence, for the above $\ve$, when $\|u_{10}-u_{20}\|$ is sufficiently small, we can guarantee that
\be\label{4.27} C_2L_{\tilde\vp}\|u_{10}-u_{20}\|\leqslant\frac\ve4.\ee
Combining the estimations from \eqref{4.20} to \eqref{4.27}, we see that for each $\ve>0$, there exists $\ol\de>0$ such that whenever $\|u_{10}-u_{20}\|<\ol\de$,
$$|\cP_t\vp(u_{10})-\cP_t\vp(u_{20})|<\ve,$$
which proves the continuity.
\eo

\subsubsection{Existence of ergodic invariant measures}

We first give the existence of invariant Borel probability measures of the original problem \eqref{1.1} -- \eqref{1.3} in $H$.
Here we use the classical procedure -- Krylov-Bogoliubov existence theorem (\cite[Corollary 11.8]{DaPZ14}) and Prokhorov's theorem (\cite[Theorem 2.3]{DaPZ14}) to show our consequence.
This procedure reads as below at our disposal.

\bl\label{le4.9} Assume that $\cP_t$ is Feller on $\cC_{\rm b}(H)$ and for every $u\in L^2(\W,H)$, define a family of Borel probability measures on $H$,
\be\label{4.28}\nu_T(\cdot)=\frac1T\int_0^T\cP_t(u,\cdot)\di t,\hs T>0.\ee
Suppose that the family $\{\nu_T\}_{T>0}$ is tight.
Then each sequence $\{\nu_{T_n}\}_{n\in\N}$ with $T_n\ra\8$ as $n\ra\8$ has a weakly convergent subsequence and the corresponding weak limit is an invariant measure for $\cP_t$.\el

For notational clarity, we denote the process given in \eqref{3.1} by $\Phi_\epsilon(t,u)$ such that $\Phi_\epsilon(t,u)(\w)=\Phi_\epsilon(t;\w,u)$ in the following discussion.
Now we show the existence of the invariant measures.

\bt\label{th4.10}Let $0<\epsilon\leqslant1$, $h\in H$ and $G$ satisfy the assumption ({\bf Gp}).
Then there exists an invariant measure in $H$ for the process $\Phi_\epsilon$.\et

\bo Following Theorem \ref{th4.8} and Lemma \ref{le4.9}, we first show the tightness of the $\nu_T$ defined as \eqref{4.28} for some suitable initial datum $u_0$, saying, given some $u_0\in L^2(\W,H)$, for arbitrary $\ve>0$, there exists a compact set $\cK_\ve\subset H$ such that
\be\label{4.29}\nu_T(\cK_\ve)\geqslant 1-\ve,\hs\mb{for all }T>0,\ee

Indeed, if we set $u_0=0$, then due to the compact embedding of $V$ into $H$ and the fact $u(t,0)\in V$ for all $t>0$, by Chebyshev's inequality and Lemma \ref{le4.3} and choosing $\cK_\ve=\ol{\bf B}_V(\cR)$, we have
\begin{align*}\nu_{T}(H\sm\ol{\bf B}_V(\cR))=&
  \nu_{T}(V\sm\ol{\bf B}_V(\cR))=\frac1{T}\int_0^{T}\bP(\|\De u(t,0)\|>\cR)\di t\\
\leqslant&
  \frac1{T\cR^2}\bE\int_0^T\|\De u(t,0)\|^2\di t\leqslant c\frac{\|h\|^2+1}{\cR^2}<\ve,
\end{align*}
as long as $\cR$ is sufficiently large.
This assures the existence of $\cK_\ve$ for \eqref{4.29}.
Therefore, Lemma \ref{le4.9} applies and the existence of invariant Borel probability measures is obtained.
\eo

At last, we cope with the existence of ergodic invariant measures for the stochastic process $\Phi_\epsilon$ generated by the problem \eqref{1.1} -- \eqref{1.3},
for which, it is necessary for us to confirm the tightness of $\cI$, that is defined to be the set of all the invariant measures for $\Phi_\epsilon$ in $H$.
We first present a boundedness estimate for the moments with respect to each invariant measure.

\bl\label{le4.11} Let $\nu\in\cI$.
Then for all $p\geqslant 1$, we have
\be\label{4.30}\int_H\|u\|^{2p}\nu(\di u)\lesssim \|h\|^{2p}+1.\ee\el
\bo We adopt the property of invariant measures to carry on the proof.
Following \eqref{4.5}, we have
\begin{align}&\di\|u(\vsig)\|^{2p}+\|u(\vsig)\|^{2p}\di\vsig+2p\|u(\vsig)\|^{2p-2}\|\De u(\vsig)\|^2\di\vsig\notag\\
=&
  -\|u\|^{2p-2}\{2p[2(\De u,u)+a\|u\|^2+(u^2,G*u^2)+(h,u)]-[p(2p-1)\epsilon^2+1]\|u\|^2\}\di\vsig\notag\\
&
  +2p\epsilon\|u\|^{2p}\di W_\vsig.\label{4.31}
\end{align}
Fix $\ol\kappa>0$ and define for each $u_0\in H$ and $p\geqslant1$ a third stopping time,
$$\ol\tau^p_{\ol{\kappa}}(u_0)=\inf\{t\geqslant0:\|u(t,u_0)\|^{2p}>\ol\kappa\}.$$

Multiplying \eqref{4.31} by $\me^\vsig$, integrating the obtained equality over $\vsig\in[0,t\wedge\ol{\tau}^p_{\ol{\kappa}}(u_0)]$ and
according to the estimates above, we deduce that
\begin{align}&\|u(t\wedge\ol{\tau}^p_{\ol{\kappa}}(u_0))\|^{2p}\notag\\
\leqslant&
  \|u_0\|^{2p}\me^{-t\wedge\ol{\tau}^p_{\ol{\kappa}}(u_0)}+\int_0^{t\wedge\ol{\tau}^p_{\ol{\kappa}}(u_0)}\me^{\vsig-t\wedge\ol{\tau}^p_{\ol{\kappa}}(u_0)}
  [\|h\|^{2p}+c\|u\|^{2p}-\al\|u\|^{2p+2}]\di\vsig\notag\\
&
  +2p\epsilon\int_0^{t\wedge\ol{\tau}^p_{\ol{\kappa}}(u_0)}\me^{\vsig-t\wedge\ol{\tau}^p_{\ol{\kappa}}(u_0)}\|u\|^{2p}\di W_\vsig.\label{4.32}
\end{align}
Let $\ol{\kappa}_0>0$ and $\|u_0\|^{2p}\leqslant\ol\kappa_0$.
Taking expectation of \eqref{4.32}, which vanishes the stochastic integral, and recalling that $c\|u\|^{2p}-\al\|u\|^{2p+2}$ has an upper bound that is independent of $u$, we obtain that
\begin{align}\bE\|u(t\wedge\ol{\tau}^p_{\ol{\kappa}}(u_0))\|^{2p}\1_{\{\|u_0\|^{2p}\leqslant\ol{\kappa}_0\}}\leqslant&
  \bE\|u_0\|^{2p}\me^{-t\wedge\ol{\tau}^p_{\ol{\kappa}}(u_0)}\1_{\{\|u_0\|^{2p}\leqslant\ol{\kappa}_0\}}\label{4.33}\\
&+c\(\|h\|^{2p}+1\)
  \bE\int_0^{t\wedge\ol{\tau}^p_{\ol{\kappa}}(u_0)}\me^{\vsig-t\wedge\ol{\tau}^p_{\ol{\kappa}}(u_0)}\1_{\{\|u_0\|^{2p}\leqslant\ol{\kappa}_0\}}\di\vsig.\notag
\end{align}
Observe that, given a fixed $\ol{\kappa}_0$, the integrands on the right-hand side of \eqref{4.33} are uniformly bounded with respect to $\kappa$, and the one on the left-hand side is non-negative.
Thus, we apply the dominated convergence theorem to the right-hand side and Fatou's lemma to the left, so that
$$\bE\|u(t)\|^{2p}\1_{\{\|u_0\|^{2p}\leqslant\ol{\kappa}_0\}}\leqslant\me^{-t}\bE\|u_0\|^{2p}\1_{\{\|u_0\|^{2p}\leqslant\ol{\kappa}_0\}}+c\(\|h\|^{2p}+1\).$$
Then for all $\ol{\kappa}_0>0$, we can pick $t_{\ol{\kappa}_0}>0$ large enough such that
\be\label{4.34}\bE\|u(t_{\ol{\kappa}_0})\|^{2p}\1_{\{\|u_0\|^{2p}\leqslant\ol{\kappa}_0\}}\leqslant2c(\|h\|^{2p}+1).\ee
Note by the property of invariant measures that, the distribution of $\|u(t_{ol{\kappa}_0})\|^{2p}$ is $\nu$.
This indicates that there exists a function $\phi_{\ol{\kappa}_0}:H\ra[0,1]$ with $\phi_{\ol{\kappa}_0}(u)\ra1$ for $\nu$-a.s. as $\ol{\kappa}_0\ra\8$, such that
\be\label{4.35}\bE\|u(t_{\ol{\kappa}_0})\|^{2p}\1_{\{\|u_0\|^{2p}\leqslant\ol{\kappa}_0\}}=\int_H\|u\|^{2p}\phi_{\ol{\kappa}_0}(u)\nu(\di u).\ee
Combining \eqref{4.34} and \eqref{4.35} and letting $\ol{\kappa}_0$ tend to the infinity, we obtain \eqref{4.30}.
\eo

\bl\label{le4.12} The invariant measure set $\cI$ is tight.\el

\bo We again follow the proof of tightness presented above.
Similarly, for every $\ve>0$, we look for a compact subset $\cK_\ve=\ol{\bf B}_V(\cR)$ with a sufficiently large $\cR$
such that
\be\label{4.36}\nu(H\sm\ol{\bf B}_V(\cR))<\ve\hs\mb{for all }\nu\in\cI.\ee
Indeed, for some fixed $t>0$, by the invariance and Lemma \ref{le4.6}, for every $\kappa>0$, we have
\begin{align}\nu(H\sm\ol{\bf B}_V(\cR))=&\int_H\cP_t(u,H\sm\ol{\bf B}_V(\cR))\nu(\di u)=\int_H\bP(\|\De\Phi_\epsilon(t,u)\|>\cR)\nu(\di u)\notag\\
=&
  \int_H\bP\(t\|\De\Phi_\epsilon(t,u)\|^2>t\cR^2\)\nu(\di u)\leqslant\int_H\bP\(\hat\tau_{t\cR^2}(u)<t\)\nu(\di u)\notag\\
\leqslant&
  \frac{C_1}{\cR^2}\left[\frac1t\int_H\|u\|^2\nu(\di u)+(\|h\|^2+1)+\|h\|^2t\right]\notag\\
&
  +\frac{ct}\kappa\left[\int_H\|u\|^4\nu(\di u)+(\|h\|^4+1)t\right],
\label{4.37}
\end{align}
where $C_1=C_1(t,\kappa)$ is given by Lemma \ref{le4.5}.
By Lemma \ref{le4.11}, we know both
$$\int_H\|u\|^2\nu(\di u)\hs\mb{and}\hs\int_H\|u\|^4\nu(\di u)\mb{ are bounded.}$$
Therefore, for each $\ve>0$, with $t$ fixed, by choosing $\kappa$ sufficiently large and fixed, we can make
$$\frac{ct}\kappa\left[\int_H\|u\|^4\nu(\di u)+(\|h\|^4+1)t\right]<\frac{\ve}2.$$
Correspondingly $C_1$ is thus fixed, and hence we can let $\cR$ be large enough such that
$$\frac{C}{\cR^2}\left[\frac1t\int_H\|u\|^2\nu(\di u)+(\|h\|^2+1)+\|h\|^2t\right]<\frac{\ve}2,$$
for the $\ve$ given above.
By the choices above, we attain an $\cR>0$ in \eqref{4.37} such that \eqref{4.36} holds.
The proof is complete.
\eo

Now we are prepared to verify the existence of ergodic invariant measures.
\bt Let $0<\epsilon\leqslant1$, $h\in H$ and $G$ satisfy the assumption ({\bf Gp}).
Then there exists an ergodic invariant measure in $H$ for $\Phi_\epsilon$.\et
\bo
By \cite[Proposition 11.12]{DaPZ14}, we see that an invariant measure is ergodic if and only if it is an extreme point of $\cI$.
We apply the Krein-Milman Theorem (see e.g. \cite[Theorem 7.4, Chapter V]{Conway90}) to show the existence of extreme points of $\cI$, i.e., to show that $\cI$ is nonempty, compact and convex.

By the discussion above, it is obvious that $\cI\ne\es$.
Note that $\cP^*_t$ is linear with respect to $\nu$, which implies that $\cI$ is convex.

For the compactness, we first show that $\cI$ is closed.
Indeed, given each weakly convergent sequence $\{\nu_n\}_{n\in\N}$ in $\cI$ such that $\nu_n$ weakly converges to $\nu$, as $n\ra\8$, by the invariance, we have for all $\Gam\in\cB(H)$,
$$
\ba{clll}\cP^*_t\nu_n(\Gam)&=\disp\int_H\cP_t(u,\Gam)\nu_n(\di u)&\ra\disp\int_H\cP_t(u,\Gam)\nu(\di u)&=\cP^*_t\nu(\Gam)\\
\|&&&\\
\nu_n(\Gam)&\ra\nu(\Gam),&\mb{as }n\ra\8.&\ea
$$
By the uniqueness of weak limit, we know that $\nu$ is also an invariant measure.
Hence the closedness.

The compactness of $\cI$ is well assured by the tightness given by Lemma \ref{le4.12} and the closedness showed above.
The proof is eventually finished.
\eo
\subsection{The case of non-negative kernel}

Now we go on to consider the existence of ergodic invariant measures of nonlocal stochastic SHE \eqref{1.1} with the special non-negative kernel, in which case $G$ satisfies \eqref{2.12}.

Actually, the only differences between the two cases lie in the estimations \eqref{4.6} and \eqref{4.32}.
The inequality \eqref{4.6} ought to be replaced in this case as follows,
\begin{align*}&\|u(s)\|^{2p}+\int_0^s\|u\|^{2p-2}\|\De u\|^2\di\vsig\\
\leqslant
&\|u_0\|^{2p}+\int_0^s[c(\|u\|^{2p}+\|h\|\|u\|^{2p-1})+\de\|u\|^{2p}-\al c_0\|u\|^{2p+2}]\di\vsig+2p\epsilon\int_0^s\|u\|^{2p}\di W_\vsig\\
\leqslant&\|u_0\|^{2p}+\int_0^s[\|h\|^{2p}+c\|u\|^{2p}-\al c_0\|u\|^{2p+2}]\di\vsig+2p\epsilon\int_0^s\|u\|^{2p}\di W_\vsig.
\end{align*}
where we have used \eqref{3.45}.
Following the same procedure of the proof of Lemma \ref{le4.3}, we can similarly obtain a conclusion of the same form as \eqref{4.4} under the assumption ({\bf Gn}).

The inequality \eqref{4.32} needs to be correspondingly modified as
\begin{align*}&\|u(t\wedge\ol{\tau}^p_{\ol{\kappa}}(u_0))\|^{2p}\\
\leqslant&
  \|u_0\|^{2p}\me^{-t\wedge\ol{\tau}^p_{\ol{\kappa}}(u_0)}+\int_0^{t\wedge\ol{\tau}^p_{\ol{\kappa}}(u_0)}\me^{\vsig-t\wedge\ol{\tau}^p_{\ol{\kappa}}(u_0)}
  [\|h\|^{2p}+c\|u\|^{2p}-\al c_0\|u\|^{2p+2}]\di\vsig\\
&
  +2p\epsilon\int_0^{t\wedge\ol{\tau}^p_{\ol{\kappa}}(u_0)}\me^{\vsig-t\wedge\ol{\tau}^p_{\ol{\kappa}}(u_0)}\|u\|^{2p}\di W_\vsig,
\end{align*}
similarly by \eqref{3.45}.
Then this is sufficient to prove a result of the same form as Lemma \ref{le4.11} if one follows the same routine in that proof.

With the same argument as that in last subsection for the rest part of the procedure, we can similarly obtain the desired conclusion as follows.
\bt Let $0<\epsilon\leqslant1$, $h\in H$ and $G$ satisfy the assumption ({\bf Gn}).
Then there exists an ergodic invariant measure in $H$ for $\Phi_\epsilon$.
\et

\section{Conclusions and remarks}

In this paper, we first show the upper semi-continuity of the $\cD$-random attractor $\sA_\epsilon(\w)$
in the regular Sobolev space $H_0^2(U)$ for the random dynamical system $\Phi_\epsilon$ generated by 2D
nonlocal stochastic Swift-Hohenberg model \eqref{1.1} -- \eqref{1.3} with multiplicative noise,
for $\bP$-a.s., $\w\in\W$, as $\epsilon\ra0^+$.
Then we obtain the existence of ergodic invariant Borel probability measures in $L^2(U)$ for the process $\Phi_\epsilon$.
What is more, both conclusions hold true for two sorts of the kernel,
i.e., the positive kernel (under assumption ({\bf Gp})) and the special non-negative kernel (under assumption ({\bf Gn})).

Firstly, we transform the original stochastic partial differential equation \eqref{1.1}
to a random partial differential equation \eqref{2.7} via a solution to an Ornstein-Uhlenbeck equation.
Then we apply the uniform estimates of the solution $v(t;\w,v_0)$ of \eqref{2.7} -- \eqref{2.9} in $H_0^2(U)$
and $\sD(A^\mu)$ for $\mu\in(\frac12,1)$ to prove the upper semi-continuity we want.
For the estimates of $v(t;\w,v_0)$ in $\sD(A^{\mu})$, we take full advantage of the analytic semigroup theory to treat the fractional power $A^{\mu}$.

Secondly, we follow the Krylov-Bogoliubov procedure to show the existence of invariant measures on $L^2(U)$,
and use \cite[Proposition 11.12]{DaPZ14} and the Krein-Milman Theorem to prove the existence of ergodic invariant measures.
Before these work, we need to assure the Feller property of the Markovian transition semigroup $\cP_t$, for which,
we apply the classical ``stochastic Gronwall's lemma", Proposition \ref{p4.1}.
However, since $h\in L^2(\W)$, when we estimate the expectation of $H_0^2(U)$-norm, one can not apply $-\De$ to \eqref{1.1} directly
and this prevents us using the It\^o's formula appropriately.
We are obliged to utilize the integral equation and the analytic semigroup theory to achieve our goal.
In this way, the classical ``stochastic Gronwall's lemma" can not directly apply at our disposal.
We hence in Lemma \ref{le4.2} develop the condition of Proposition \ref{p4.1} into an equivalent condition that is much easier to check.
This enables us to handle this difficulty.

Thirdly, as an extension, our approach of proving the existence of ergodic invariant measures can also apply to more general cases.
For example, the stochastic term $\epsilon u\di W_t$ in \eqref{1.1} can be replaced by the formal expansion
$$\psi(u)\di\~W_t=\sum_{k\geqslant1}\psi_k(u)\di W^k_t,$$
where $\{W_t^k\}_{k\geqslant1}$ is a sequence of independent one-dimensional Wiener process relative to some prescribed stochastic basis and
$\psi$ is required to be naturally compatible, bounded and Lipschitz continuous.
Also, our technique of the application of analytic semigroup theory enlightens us to deal with more general stochastic partial differential equations, such as non-autonomous ones.

At last, we would like to raise some interesting problems related to this paper.
\benu\item[(1)]We have seen the existences of $\cD$-random attractors $\sA_\epsilon(\w)$, for $\bP$-a.s. $\w\in\W$ and the invariant measures for the process $\Phi_\epsilon$.
We know that for deterministic dynamical systems, the invariant measures are often supported on (global, pullback) attractors.
Then what is the relationship between the random attractors and invariant measures for stochastic dynamical system $\Phi_\epsilon$?
\item[(2)]Now that we have obtained the existence of invariant measures on $L^2(U)$,
the next challenging work is to improve the regularity of the invariant measures for the nonlocal stochastic Swift-Hohenberg model \eqref{1.1} -- \eqref{1.3}.
And it is also of great interest to determine the uniqueness and highest regularity of the invariant measures.
\item[(3)]In comparison with the upper semi-continuity of $\cD$-random attractor $\sA_\epsilon(\w)$ as $\epsilon\ra0^+$,
do the sets $\cI_\epsilon$ of all invariant measures for the process $\Phi_\epsilon$ have upper semi-continuity as $\epsilon\ra0^+$?
If so, whether are they upper semi-continuous to the set of all invariant measures of the deterministic dynamical systems given by \eqref{3.47}?
\eenu

\section*{Data Availability}

Data sharing is not applicable to this article as no new data were created or analyzed in this study.

\section*{Acknowledgements}

Our work was supported by grants from the National Natural Science Foundation of China (NSFC No. 11801190 and 12071192).


\end{document}